\documentclass[a4paper]{article}

\usepackage{amssymb}

\def\GAP{\textsf{GAP}}
\def\ATLAS{\textsc{Atlas}}
\def\MAGMA{\textsf{MAGMA}}

\def\Cent{{\rm\bf C}}
 \def\Z{{\mathbb Z}} 
  
\def\B{{\mathbb B}}
\def\tthdump#1{#1}
\tthdump{}

\begin{document}

\tthdump{\title{Some steps in the verification of the ordinary character table of the Baby Monster group}}

\author{\textsc{Thomas Breuer, Kay Magaard, Robert A.~Wilson}}

\date{May 17th, 2019}

\maketitle

\abstract{%
We show the details of certain computations that are described
in \cite{BMverify}.}

\textwidth16cm
\oddsidemargin0pt

\parskip 1ex plus 0.5ex minus 0.5ex
\parindent0pt

\tableofcontents


\section{Overview}

The aim of~\cite{BMverify} is to verify the ordinary character table
of the Baby Monster group $\B$.
Here we collect,
in the form of an explicit {\GAP}~\cite{GAP} session protocol,
the computations that are needed in that paper.

We proceed as follows.

Section~\ref{pres_B}
shows the computations that are described in~\cite[Section~3]{BMverify}.
At this point, we know that the three matrix groups that are used
later on are in fact representations of the group $\B$,
w.~r.~t.~compatible (standard) generators.

Section~\ref{invs_B} turns the class invariants and the power map
information from~\cite[Section~4]{BMverify} into a {\GAP} function
that identifies the class label of a given word in terms of the
given standard generators.

Section~\ref{centralizers_prime}
shows part of the computations described in~\cite[Section~5]{BMverify}.

Section~\ref{table_c2b} shows the computation of the character table
of an involution centralizer of type $2^{1+22}.Co_2$ in $\B$.

Section~\ref{sect:classes} shows how the conjugacy classes,
the corresponding centralizer orders, and the power maps of $\B$
are determined.

In Section~\ref{sect:classes}, we put these pieces together and write down
the list of class representatives of $\B$,
together with their centralizer orders and power maps.

With this information and with the (already verified) character tables
of some known subgroups of $\B$,
computing the irreducible characters of $\B$ is then easy;
this corresponds to~\cite[Section~7]{BMverify},
and is done in Section~\ref{sect:irreducibles}.

We will use the {\GAP} Character Table Library
and the interface to the {\ATLAS} of Group Representations~\cite{AGRv3},
thus we load these {\GAP} packages.

\begin{verbatim}
    gap> LoadPackage( "ctbllib", false );
    true
    gap> LoadPackage( "atlasrep", false );
    true
\end{verbatim}

The {\MAGMA}~\cite{Magma} system will be needed
for computing a character table
and for several conjugacy tests.
If the following command returns \texttt{false}
then these steps will not work.

\begin{verbatim}
    gap> CTblLib.IsMagmaAvailable();
    true
\end{verbatim}

\section{Verification of a presentation for $\B$}\label{pres_B}

We show the computations that are described in~\cite[Section~4]{BMverify}.
First we create the free generators and relators of the presentation.

\begin{verbatim}
    gap> F:= FreeGroup( List( "abcdefghijk", x -> [ x ] ) );;
    gap> gens:= GeneratorsOfGroup( F );;
    gap> rels:= List( gens, x -> x^2 );;
    gap> ord3pairs:= List( [ 1 .. 7 ], i -> [ i, i+1 ] );
    [ [ 1, 2 ], [ 2, 3 ], [ 3, 4 ], [ 4, 5 ], [ 5, 6 ], [ 6, 7 ], [ 7, 8 ] ]
    gap> Append( ord3pairs, [ [ 5, 9 ], [ 9, 10 ], [ 10, 11 ] ] );
    gap> for pair in ord3pairs do
    >      Add( rels, ( gens[ pair[1] ] * gens[ pair[2] ] )^3 );
    >    od;
    gap> for i in [ 1 .. 11 ] do
    >      for j in [ i+1 .. 11 ] do
    >        if not [ i, j ] in ord3pairs then
    >          Add( rels, ( gens[i] * gens[j] )^2 );
    >        fi;
    >      od;
    >    od;
    gap> Add( rels, Product( gens{ [ 5, 4, 3, 5, 6, 7, 5, 9, 10 ] } )^10 );
\end{verbatim}

We do not call \texttt{FreeGroup( 11 )} because later on we want to
translate the relators into straight line programs,
and we can use \texttt{StraightLineProgram} with first argument a string
only if no generator name is a prefix of another generator name.

\begin{verbatim}
    gap> gensstrings:= List( gens, String );;
    gap> relsslps:= List( List( rels, String ),
    >                     x -> StraightLineProgram( x, gensstrings ) );;
\end{verbatim}

Next we write a straight line program that computes the $11$ generators
$t_1, \ldots, t_{11}$,
following the steps shown in~\cite[Table~1]{BMverify}.
We start with the two standard generators $a$ and $b$, say,
in the slots 1 and 2,
and compute expressions for the subsequent slots.
The product $a b$ will be in position 3, its 5th power $(a b)^5$
(which will be needed later on) in position 4,
the power $(a b)^{15}$ in position 5, and $d = (a b)^{15} b$
in position 6.
The generators $t_{11} = d^{19}$ gets stored in position 7.

\begin{verbatim}
    gap> slp:= [ [ 1, 1, 2, 1 ], [ 3, 5 ], [ 4, 3 ], [ 5, 1, 2, 1 ] ];;
    gap> resultpos:= [];;
    gap> Add( slp, [ 6, 19 ] );
    gap> resultpos[11]:= Length( slp ) + 2;;
\end{verbatim}

Next we compute $c = (a t_{11})^3$ (position 9),
$e = ((c d^3)^{10})^d$ (position 13), \ldots

\begin{verbatim}
    gap> Append( slp, [ [ 1, 1, 7, 1 ], [ 8, 3 ] ] );
    gap> Append( slp, [ [ 6, 3 ], [ 9, 1, 10, 1 ], [ 11, 10 ],
    >                   [ 6, -1, 12, 1, 6, 1 ] ] );
\end{verbatim}

\ldots $t_1 = f = ((((e c)^6 c (e c)^3 )^2 e c e^2 c)^5)^((e c)^4)$
(position 24), \ldots

\begin{verbatim}
    gap> # 14: e*c,  15: (e*c)^2,  16: (e*c)^3,  17: (e*c)^4,  18: (e*c)^6
    gap> Append( slp, [ [ 13, 1, 9, 1 ], [ 14, 2 ], [ 14, 1, 15, 1 ],
    >                   [ 15, 2 ], [ 16, 2 ] ] );
    gap> # 19: e*c*e,  20: e*c*e^2*c
    gap> Append( slp, [ [ 14, 1, 13, 1 ], [ 19, 1, 14, 1 ] ] );
    gap> # 21: (e*c)^6*c*(e*c)^3,  22: ((e*c)^6*c*(e*c)^3)^2*e*c*e^2*c
    gap> Append( slp, [ [ 18, 1, 9, 1, 16, 1 ], [ 21, 2, 20, 1 ] ] );
    gap> # 23: (((e*c)^6*c*(e*c)^3)^2*e*c*e^2*c)^5
    gap> Append( slp, [ [ 22, 5 ] ] );
    gap> # 24: t1 = f = ((((e*c)^6*c*(e*c)^3)^2*e*c*e^2*c)^5)^((e*c)^4)
    gap> Append( slp, [ [ 17, -1, 23, 1, 17, 1 ] ] );
    gap> resultpos[1]:= Length( slp ) + 2;;
\end{verbatim}

\ldots $g = ((e c)^8 c (e c)^3 )^{{(e c e^2 c)^2}}$
(position 27) and $t_2 = f^{{g^f}}$ (position 30), \ldots

\begin{verbatim}
    gap> # 25: (e*c)^8*c*(e*c)^3,  26: (e*c*e^2*c)^2
    gap> Append( slp, [ [ 15, 1, 21, 1 ], [ 20, 2 ] ] );
    gap> # 27: g = ((e*c)^8*c*(e*c)^3)^((e*c*e^2*c)^2)
    gap> Append( slp, [ [ 26, -1, 25, 1, 26, 1 ] ] );
    gap> # 28: g f,  29: g^-1,  30: t2 = f^{{g f}}
    gap> Append( slp, [ [ 27, 1, 24, 1 ], [ 27, -1 ],
    >                   [ 28, -1, 24, 1, 28, 1 ] ] );
    gap> resultpos[2]:= Length( slp ) + 2;;
\end{verbatim}

\ldots $t_3 = f^{{g f g}}$,
$t_4 = f^{{g f g^2}}$,
$t_5 = f^{{g f g^3}}$,
$t_6 = f^{{g f g^4}}$,
$t_7 = f^{{g f g^5}}$,
$t_8 = f^{{g f g^6}}$ (positions 31 to 36), \ldots

\begin{verbatim}
    gap> # 31: t3 = f^( g * f * g )
    gap> Append( slp, [ [ 29, 1, 30, 1, 27, 1 ] ] );
    gap> resultpos[3]:= Length( slp ) + 2;;
    gap> # 32: t4 = f^( g * f * g^2 )
    gap> Append( slp, [ [ 29, 1, 31, 1, 27, 1 ] ] );
    gap> resultpos[4]:= Length( slp ) + 2;;
    gap> # 33: t5 = f^( g * f * g^3 )
    gap> Append( slp, [ [ 29, 1, 32, 1, 27, 1 ] ] );
    gap> resultpos[5]:= Length( slp ) + 2;;
    gap> # 34: t6 = f^( g * f * g^4 )
    gap> Append( slp, [ [ 29, 1, 33, 1, 27, 1 ] ] );
    gap> resultpos[6]:= Length( slp ) + 2;;
    gap> # 35: t7 = f^( g * f * g^5 )
    gap> Append( slp, [ [ 29, 1, 34, 1, 27, 1 ] ] );
    gap> resultpos[7]:= Length( slp ) + 2;;
    gap> # 36: t8 = f^( g * f * g^6 )
    gap> Append( slp, [ [ 29, 1, 35, 1, 27, 1 ] ] );
    gap> resultpos[8]:= Length( slp ) + 2;;
\end{verbatim}

\ldots $p = ((a b)^5 t_{11} (a b)^{-5} t_1 (a b)^5)^{-1}$ (position 38)
and $i = d^p$ (position 41), \ldots

\begin{verbatim}
    gap> # 37: (a*b)^5*t11*(a*b)^-5*t1*(a*b)^5,
    gap> # 38: p = ((a*b)^5*t11*(a*b)^-5*t1*(a*b)^5)^-1
    gap> Append( slp, [ [ 4, 1, 7, 1, 4, -1, 24, 1, 4, 1 ], [ 37, -1 ] ] );
    gap> # 39: p^-1,  40: h = c^p, 41: i = d^p
    gap> Append( slp, [ [ 38, -1 ], [ 39, 1, 9, 1, 38, 1 ],
    >                   [ 39, 1, 6, 1, 38, 1 ] ] );
\end{verbatim}

\ldots $j = [t_5^{{i^2}}, t_3 t_4]$ (position 45) and
$k = [t_5^{{i^5}}, t_3 t_4]$ (position 48), \ldots

\begin{verbatim}
    gap> # 42: i^2,  43: t5^(i^2),  44: t3*t4
    gap> Append( slp, [ [ 41, 2 ], [ 42, -1, 33, 1, 42, 1 ],
    >                   [ 31, 1, 32, 1 ] ] );
    gap> # 45: j = Comm( t5^( i^2 ), t3*t4 )
    gap> Append( slp, [ [ 43, -1, 44, -1, 43, 1, 44, 1 ] ] );
    gap> # 46: i^3,  47: t5^(i^5),  48: k = Comm( t5^(i^5), t3*t4 )
    gap> Append( slp, [ [ 41, 1, 42, 1 ], [ 46, -1, 43, 1, 46, 1 ],
    >                   [ 47, -1, 44, -1, 47, 1, 44, 1 ] ] );
\end{verbatim}

\ldots $l = [t_8^{{j k}}, t_6 t_7] [t_8^{{k j}}, t_6 t_7]$ (position 57),
\ldots

\begin{verbatim}
    gap> # 49: t6*t7,  50: (t6*t7)^-1,  51: j*k,  52: k*j,  53: t8^(j*k)
    gap> Append( slp, [ [ 34, 1, 35, 1 ], [ 49, -1 ], [ 45, 1, 48, 1 ],
    >                   [ 48, 1, 45, 1 ], [ 51, -1, 36, 1, 51, 1 ] ] );
    gap> # 54: Comm( t8^(j*k), t6*t7),  55: t8^(k*j)
    gap> Append( slp, [ [ 53, -1, 50, 1, 53, 1, 49, 1 ] ] );
    gap> Append( slp, [ [ 52, -1, 36, 1, 52, 1 ] ] );
    gap> # 56: Comm( t8^(k*j), t6*t7 )
    gap> Append( slp, [ [ 55, -1, 50, 1, 55, 1, 49, 1 ] ] );
    gap> # 57: l = Comm( t8^(j*k), t6*t7 ) * Comm( t8^(k*j), t6*t7 )
    gap> Append( slp, [ [ 54, 1, 56, 1 ] ] );
\end{verbatim}

\ldots $l_3 = [t_8^{{(j k)^4}}, t_6 t_7]$ (position 61),
$l_4 = t_8^{{(j k)^3 k j}}$ (position 62),
$l_5 = (l l_3 l_4)^3 l_3 l_4$ (position 65), \ldots

\begin{verbatim}
    gap> # 58: (j*k)^3,  59: (j*k)^-3,  60: t8^((j*k)^4)
    gap> Append( slp, [ [ 51, 3 ], [ 58, -1 ], [ 59, 1, 53, 1, 58, 1 ] ] );
    gap> # 61: l3 = Comm( t8^((j*k)^4), t6*t7 )
    gap> Append( slp, [ [ 60, -1, 50, 1, 60, 1, 49, 1 ] ] );
    gap> # 62: l4 = t8^((j*k)^3*k*j)
    gap> Append( slp, [ [ 52, -1, 59, 1, 36, 1, 58, 1, 52, 1 ] ] );
    gap> # 63: l3*l4,  64: l*l3*l4
    gap> Append( slp, [ [ 61, 1, 62, 1 ], [ 57, 1, 63, 1 ] ] );
    gap> # 65: l5:= ( l * l3 * l4 )^3 * l3 * l4;;
    gap> Append( slp, [ [ 64, 3, 63, 1 ] ] );
\end{verbatim}

\ldots $m_2 = l_4^{{l_5^4}}$ (position 67),
$m_3 = m_2^{{l_5}}$ (position 68),
$t_{10} = m_3 m_2 l_4 m_2 m_3$ (position 69), and
$t_9 = l_4 m_2 t_{10} m_2 l_4$ (position 70).

\begin{verbatim}
    gap> # 66: l5^4,  67: m2 = l4^(l5^4)
    gap> Append( slp, [ [ 65, 4 ], [ 66, -1, 62, 1, 66, 1 ] ] );
    gap> # 68: m3 = m2^l5
    gap> Append( slp, [ [ 65, -1, 67, 1, 65, 1 ] ] );
    gap> # 69: t10 = m3*m2*l4*m2*m3
    gap> Append( slp, [ [ 68, 1, 67, 1, 62, 1, 67, 1, 68, 1 ] ] );
    gap> resultpos[10]:= Length( slp ) + 2;;
    gap> # 70: t9 = l4*m2*t10*m2*l4
    gap> Append( slp, [ [ 62, 1, 67, 1, 69, 1, 67, 1, 62, 1 ] ] );
    gap> resultpos[9]:= Length( slp ) + 2;;
\end{verbatim}

Finally, we specify the list of outputs,
and create the straight line program object.

\begin{verbatim}
    gap> Add( slp, List( resultpos, x -> [ x, 1 ] ) );
    gap> slp:= StraightLineProgram( slp, 2 );
    <straight line program>
\end{verbatim}

And now we compute,
for each of the three pairs of generators we are interested in,
the $11$ generators,
and test whether these generators satisfy the presentation.

\begin{verbatim}
    gap> b_2:= AtlasGroup( "B", Characteristic, 2, Dimension, 4370 );;
    gap> b_3:= AtlasGroup( "B", Characteristic, 3, Dimension, 4371 );;
    gap> b_5:= AtlasGroup( "B", Characteristic, 5, Dimension, 4371 );;
    gap> gens_2:= GeneratorsOfGroup( b_2 );;
    gap> gens_3:= GeneratorsOfGroup( b_3 );;
    gap> gens_5:= GeneratorsOfGroup( b_5 );;
    gap> res_2:= ResultOfStraightLineProgram( slp, gens_2 );;
    gap> ForAll( relsslps,
    >            prg -> IsOne( ResultOfStraightLineProgram( prg, res_2 ) ) );
    true
    gap> res_3:= ResultOfStraightLineProgram( slp, gens_3 );;
    gap> ForAll( relsslps,
    >            prg -> IsOne( ResultOfStraightLineProgram( prg, res_3 ) ) );
    true
    gap> res_5:= ResultOfStraightLineProgram( slp, gens_5 );;
    gap> ForAll( relsslps,
    >            prg -> IsOne( ResultOfStraightLineProgram( prg, res_5 ) ) );
    true
\end{verbatim}

In order to prove that the $11$ elements that satisfy the relations
generate the same group as the original generators,
we create a straight line program that computes the elements $a', b'$
stated in \cite[Section 4.4]{BMverify},
first the elements $r$ and $s$ (positions 12 and 13), \ldots

\begin{verbatim}
    gap> revslp:= [ Concatenation( List( [ 1 .. 8 ], i -> [ i, 1 ] ) ),
    >        Concatenation( List( [ 5, 9, 10, 11 ], i -> [ i, 1 ] ) ) ];;
\end{verbatim}

\ldots and then $a' = (r^7 s)^{15}$ (position 15)
and $b' = (t_1 t_2)^{(sr)^10}$ (position 20).

\begin{verbatim}
    gap> Append( revslp, [ [ 12, 7, 13, 1 ], [ 14, 15 ], [ 13, 1, 12, 1 ],
    >                      [ 16, 10 ], [ 17, -1 ], [ 1, 1, 2, 1 ],
    >                      [ 18, 1, 19, 1, 17, 1 ] ] );
\end{verbatim}

Again, we specify the outputs and create the straight line program object.

\begin{verbatim}
    gap> Add( revslp, List( [ 15, 20 ], x -> [ x, 1 ] ) );
    gap> revslp:= StraightLineProgram( revslp, 11 );
    <straight line program>
\end{verbatim}

We claim that, for the three representations in question,
evaluating the straight line program \texttt{revslp} at the $11$
generators
yields a pair $a', b'$ of matrices that is simultaneously conjugate to the
original matrices $a, b$.
Once this is established, we know that the group $\langle a, b \rangle$
is equal to the group generated by the $11$ generators,
and that mapping the original generators of any of the three
representations to the original generators of another one
defines a group isomorphism.

In order to show the conjugacy property,
we use that the nullspace of $w(a, b) = a^0 + a b + b a + b$
is $1$-dimensional, in all three cases.

\begin{verbatim}
    gap> a:= gens_2[1];; b:= gens_2[2];;
    gap> w:= One( a ) + b*a + a*b + b;;
    gap> nsp_2:= NullspaceMat( w );; Length( nsp_2 );
    1
    gap> a:= gens_3[1];; b:= gens_3[2];;
    gap> w:= One( a ) + b*a + a*b + b;;
    gap> nsp_3:= NullspaceMat( w );; Length( nsp_3 );
    1
    gap> a:= gens_5[1];; b:= gens_5[2];;
    gap> w:= One( a ) + b*a + a*b + b;;
    gap> nsp_5:= NullspaceMat( w );; Length( nsp_5 );
    1
\end{verbatim}

The standard basis w.~r.~t.~given generators and a vector $v$
is defined by starting with the list $b = [ v ]$
and iteratively adding those images of the vectors in $b$
under the right multiplication with the generators that increase the
dimension of the vector space generated by $b$.
(Since such a function is apparently not available in {\GAP}'s MeatAxe,
we provide it here.)

\begin{verbatim}
    gap> StdBasis:= function( F, mats, seed )
    >      local n, b, mb, v, m, new;
    > 
    >      n:= Length( mats[1] );
    >      b:= [ seed ];
    >      mb:= MutableBasis( F, b );
    >      for v in b do
    >        for m in mats do
    >          new:= v * m;
    >          if not IsContainedInSpan( mb, new ) then
    >            Add( b, new );
    >            if Length( b ) = n then
    >              break;
    >            fi;
    >            CloseMutableBasis( mb, new );
    >          fi;
    >        od;
    >        if Length( b ) = n then
    >          break;
    >        fi;
    >      od;
    >      return b;
    >    end;;
\end{verbatim}


All we have to check is that the matrices of the linear mappings $a, b$
w.~r.~t.~their standard basis and a generating vector of the nullspace
of $w(a, b)$ are equal to the matrices of $a', b'$ w.~r.~t.~their
standard basis and a generating vector of the nullspace of $w(a', b')$.

We verify this in characteristic $2$, \ldots

\begin{verbatim}
    gap> stdbas_2:= StdBasis( GF(2), gens_2, nsp_2[1] );;
    gap> inv:= stdbas_2^-1;;
    gap> stdgens_2:= List( gens_2, m -> stdbas_2 * m * inv );;
    gap> newgens_2:= ResultOfStraightLineProgram( revslp, res_2 );;
    gap> aa:= newgens_2[1];;  bb:= newgens_2[2];;
    gap> neww:= One( aa ) + bb * aa + aa * bb + bb;;
    gap> newnsp_2:= NullspaceMat( neww );;  Length( newnsp_2 );
    1
    gap> newstdbas_2:= StdBasis( GF(2), newgens_2, newnsp_2[1] );;
    gap> inv:= newstdbas_2^-1;;
    gap> newstdgens_2:= List( newgens_2, m -> newstdbas_2 * m * inv );;
    gap> stdgens_2 = newstdgens_2;
    true
\end{verbatim}

\ldots in characteristic $3$, \ldots

\begin{verbatim}
    gap> stdbas_3:= StdBasis( GF(3), gens_3, nsp_3[1] );;
    gap> inv:= stdbas_3^-1;;
    gap> stdgens_3:= List( gens_3, m -> stdbas_3 * m * inv );;
    gap> newgens_3:= ResultOfStraightLineProgram( revslp, res_3 );;
    gap> aa:= newgens_3[1];;  bb:= newgens_3[2];;
    gap> neww:= One( aa ) + bb * aa + aa * bb + bb;;
    gap> newnsp_3:= NullspaceMat( neww );;  Length( newnsp_3 );
    1
    gap> newstdbas_3:= StdBasis( GF(3), newgens_3, newnsp_3[1] );;
    gap> inv:= newstdbas_3^-1;;
    gap> newstdgens_3:= List( newgens_3, m -> newstdbas_3 * m * inv );;
    gap> stdgens_3 = newstdgens_3;
    true
\end{verbatim}

\ldots and in characteristic $5$.

\begin{verbatim}
    gap> stdbas_5:= StdBasis( GF(5), gens_5, nsp_5[1] );;
    gap> inv:= stdbas_5^-1;;
    gap> stdgens_5:= List( gens_5, m -> stdbas_5 * m * inv );;
    gap> newgens_5:= ResultOfStraightLineProgram( revslp, res_5 );;
    gap> aa:= newgens_5[1];;  bb:= newgens_5[2];;
    gap> neww:= One( aa ) + bb * aa + aa * bb + bb;;
    gap> newnsp_5:= NullspaceMat( neww );;  Length( newnsp_5 );
    1
    gap> newstdbas_5:= StdBasis( GF(5), newgens_5, newnsp_5[1] );;
    gap> inv:= newstdbas_5^-1;;
    gap> newstdgens_5:= List( newgens_5, m -> newstdbas_5 * m * inv );;
    gap> stdgens_5 = newstdgens_5;
    true
\end{verbatim}

\section{Invariants that distinguish conjugacy classes of $\B$}\label{invs_B}

The function \texttt{IdentifyClassName} shown below
implements the invariants defined in~\cite[Section~5]{BMverify},
that distinguish $183$ conjugacy classes of $\B$.

Its input can be as follows.

\begin{itemize}
\item
  Three matrices \texttt{data2}, \texttt{data3}, \texttt{data5},
  representing an element of $\B$ in the given three matrix representations,
  in characteristics $2$, $3$, and $5$, respectively;
  in this case, the argument \texttt{slp} should be \texttt{fail}.
\item
  Lists \texttt{data2}, \texttt{data3}, \texttt{data5}
  of standard generators of $\B$ in the given three matrix representations
  such that \texttt{slp} is a straight line program
  that takes these generators as inputs,
  and computes the element in question.
\end{itemize}

In both cases, the argument \texttt{order} can be either \texttt{fail}
or the order of the element.
A known order allows us to omit any computation with matrices
in several cases.

The output is the label for the union of conjugacy classes
as defined in~\cite[Table 2]{BMverify},
except that labels containing \emph{two} letters are returned
in those cases that will later turn out to describe two Galois conjugate
classes
--these are \texttt{"23AB"}, \texttt{"30GH"}, \texttt{"31AB"},
\texttt{"32AB"}, \texttt{"32CD"}, \texttt{"34BC"}, \texttt{"46AB"},
\texttt{"47AB"}, \texttt{"56AB"}--
and in the case of \texttt{"16DF"} where we have no invariant that
distinguishes two classes that are not Galois conjugate.

\begin{verbatim}
    gap> IdentifyClassName:= function( data2, data3, data5, slp, order )
    >      local data, mats, elm, cand, nams, trace, pos, one, rank;
    > 
    >      data:= [ , data2, data3,, data5 ];
    >      mats:= [];
    > 
    >      elm:= function( p )
    >        if not IsBound( mats[p] ) then
    >          if slp = fail then
    >            mats[p]:= data[p];
    >          else
    >            mats[p]:= ResultOfStraightLineProgram( slp, data[p] );
    >          fi;
    >        fi;
    >        return mats[p];
    >      end;
    > 
    >      if order = fail then
    >        order:= Order( elm(2) );
    >      fi;
    > 
    >      # The element order suffices in certain cases.
    >      if order in [ 23, 31, 46, 47, 56 ] then
    >        # There are two Galois conjugate classes of elements of this order.
    >        return Concatenation( String( order ), "AB" );
    >      elif order in [ 1, 7, 11, 13, 17, 19, 21, 25, 27, 33, 35, 38, 39,
    >                      44, 47, 52, 55, 66, 70 ] then
    >        # There is exactly one conjugacy class of elements of this order.
    >        return Concatenation( String( order ), "A" );
    >      fi;
    > 
    >      if order in [ 3, 5, 9, 15 ] then
    >        # The trace in the 2-modular representation suffices.
    >        cand:= [ [3,1], [3,0], [5,0], [5,1], [9,0], [9,1], [15,0], [15,1] ];
    >        nams:= [ "3A", "3B", "5A", "5B", "9A", "9B", "15A", "15B" ];
    >        trace:= Int( TraceMat( elm(2) ) );
    >        return nams[ Position( cand, [ order, trace ] ) ];
    >      elif order mod 4 = 2 then
    >        # Compute the rank of 1 + x.
    >        cand:= [ [ 2, 1860 ], [ 2, 2048 ], [ 2, 2158 ], [ 2, 2168 ],
    >                 [ 6, 3486 ], [ 6, 3510 ], [ 6, 3566 ], [ 6, 3534 ],
    >                 [ 6, 3606 ], [ 6, 3604 ], [ 6, 3596 ], [ 6, 3610 ],
    >                 [ 6, 3636 ], [ 6, 3638 ], [ 6, 3634 ],
    >                 [ 10, 3860 ], [ 10, 3896 ], [ 10, 3918 ], [ 10, 3908 ],
    >                 [ 10, 3920 ], [ 10, 3932 ],
    >                 [ 14, 3996 ], [ 14, 4008 ], [ 14, 4048 ], [ 14,4034 ],
    >                 [ 14, 4052 ],
    >                 [ 18, 4088 ], [ 18, 4090 ], [ 18, 4110 ], [ 18, 4124 ],
    >                 [ 18, 4128 ], [ 18, 4122 ],
    >                 [ 22, 4140 ], [ 22, 4158 ],
    >                 [ 26, 4198 ], [ 26, 4176 ],
    >                 [ 30, 4190 ], [ 30, 4212 ], [ 30, 4206 ], [ 30, 4214 ],
    >                 [ 30, 4224 ], [ 30, 4216 ],
    >                 [ 34, 4238 ], [ 34, 4220 ],
    >                 [ 42, 4242 ], [ 42, 4258 ] ];
    >        nams:= [ "2A", "2B", "2C", "2D",
    >                 "6A", "6B", "6C", "6D", "6E", "6F",
    >                 "6G", "6H", "6I", "6J", "6K",
    >                 "10A", "10B", "10C", "10D", "10E", "10F",
    >                 "14A", "14B", "14C", "14D", "14E",
    >                 "18A", "18B", "18C", "18D", "18E", "18F",
    >                 "22A", "22B",
    >                 "26A", "26B",
    >                 "30AB", "30C", "30D", "30E", "30F", "30GH",
    >                 "34A", "34BC",
    >                 "42AB", "42C" ];
    >        one:= elm(2)^0;
    >        rank:= RankMat( elm(2) + one );
    >        pos:= Position( cand, [ order, rank ] );
    >        if nams[ pos ] = "30AB" then
    >          rank:= RankMat( elm(2)^5 + one );
    >          if   rank = 3510 then return "30A";
    >          elif rank = 3486 then return "30B";
    >          else Error( "wrong rank" );
    >          fi;
    >        elif nams[ pos ] = "42AB" then
    >          rank:= RankMat( elm(2)^3 + one );
    >          if   rank = 3996 then return "42A";
    >          elif rank = 4008 then return "42B";
    >          else Error( "wrong rank" );
    >          fi;
    >        else
    >          return nams[ pos ];
    >        fi;
    >      elif order in [ 36, 60 ] then
    >        cand:= [ [36,4226],[36,4238],[36,4248],[60,4280],[60,4286],[60,4296] ];
    >        nams:= [ "36A", "36B", "36C", "60A", "60B", "60C" ];
    >        one:= elm(2)^0;
    >        rank:= RankMat( elm(2) + one );
    >        return nams[ Position( cand, [ order, rank ] ) ];
    >      elif order = 28 then
    >        one:= elm(2)^0;
    >        rank:= RankMat( elm(2) + one );
    >        trace:= Int( TraceMat( elm(3)^7 ) );
    >        if rank = 4188 then   # 28A or 28C
    >          if   trace = 0 then return "28A";
    >          elif trace = 1 then return "28C";
    >          else Error( "wrong trace" );
    >          fi;
    >        elif rank = 4200 then   # 28B or 28D
    >          if   trace = 1 then return "28B";
    >          elif trace = 0 then return "28D";
    >          else Error( "wrong trace" );
    >          fi;
    >        elif rank = 4210 then return "28E";
    >        else Error( "wrong rank" );
    >        fi;
    >      elif order = 32 then
    >        trace:= Int( TraceMat( elm(3)^2 ) );
    >        if   trace = 2 then return "32AB";
    >        elif trace = 0 then return "32CD";
    >        else Error( "wrong trace" );
    >        fi;
    >      elif order = 40 then
    >        one:= elm(2)^0;
    >        rank:= RankMat( elm(2) + one );
    >        if rank = 4242 then   # 40A or 40B 0r 40C
    >          trace:= Int( TraceMat( elm(3) ) );
    >          if   trace = 0 then return "40A";
    >          elif trace = 1 then return "40B";
    >          else return "40C";
    >          fi;
    >        elif rank = 4250 then return "40D";
    >        elif rank = 4258 then return "40E";
    >        else Error( "wrong rank" );
    >        fi;
    >      elif order = 48 then
    >        trace:= Int( TraceMat( elm(3) ) );
    >        if   trace = 0 then return "48A";
    >        elif trace = 1 then return "48B";
    >        else Error( "wrong trace" );
    >        fi;
    >      elif order in [ 4, 8 ] then
    >        cand:= [ [4,3114],[4,3192],[4,3256],[4,3202],[4,3204],[4,3266],
    >                 [4,3264],[8,3774],[8,3738],[8,3778],[8,3780],[8,3810],
    >                 [8,3786],[8,3812],[8,3818] ];
    >        nams:= [ "4A-B", "4C-D", "4E", "4F", "4G", "4H-J", "4I",
    >                 "8A", "8B-C-E", "8D", "8F-H", "8G", "8I-L", "8J", "8K-M-N" ];
    >        one:= elm(2)^0;
    >        rank:= RankMat( elm(2) + one );
    >        pos:= Position( cand, [ order, rank ] );
    >        if not '-' in nams[ pos ] then
    >          return nams[ pos ];
    >        elif order = 4 then
    >          trace:= Int( TraceMat( elm(3) ) );
    >          if trace = 0 then
    >            if   nams[ pos ] = "4A-B" then return "4B";
    >            elif nams[ pos ] = "4C-D" then return "4D";
    >            else return "4H";
    >            fi;
    >          elif trace = 1 then
    >            if   nams[ pos ] = "4A-B" then return "4A";
    >            elif nams[ pos ] = "4C-D" then return "4C";
    >            else return "4J";
    >            fi;
    >          else
    >            Error( "wrong trace" );
    >          fi;
    >        elif order = 8 then
    >          if nams[ pos ] = "8B-C-E" then
    >            trace:= Int( TraceMat( elm(3) ) );
    >            if   trace = 1 then return "8B";
    >            elif trace = 0 then return "8C";
    >            else return "8E";
    >            fi;
    >          elif nams[ pos ] = "8F-H" then
    >            rank:= RankMat( ( elm(2) + one )^3 );
    >            if   rank = 2619 then return "8F";
    >            elif rank = 2620 then return "8H";
    >            else Error( "wrong rank" );
    >            fi;
    >          elif nams[ pos ] = "8I-L" then
    >            rank:= RankMat( ( elm(2) + one )^2 );
    >            if   rank = 3202 then return "8I";
    >            elif rank = 3204 then return "8L";
    >            else Error( "wrong rank" );
    >            fi;
    >          else   # 8K-M-N
    >            rank:= RankMat( ( elm(2) + one )^3 );
    >            if rank = 2714 then   # 8K-M
    >              trace:= Int( TraceMat( elm(3) ) );
    >              if   trace = 1 then return "8K";
    >              elif trace = 2 then return "8M";
    >              else Error( "wrong trace" );
    >              fi;
    >            elif rank = 2717 then return "8N";
    >            else Error( "wrong rank" );
    >            fi;
    >          fi;
    >        fi;
    >      elif order in [ 12, 24 ] then
    >        cand:= [ [12,3936],[12,3942],[12,3958],[12,3996],[12,3962],[12,3964],
    >                 [12,3986],[12,3978],[12,3966],[12,4000],[12,3982],[12,3988],
    >                 [12,4002],[12,4004],[24,4152],[24,4164],[24,4170],[24,4182],
    >                 [24,4176],[24,4178],[24,4174],[24,4186] ];
    >        nams:= [ "12A-C-D", "12B", "12E", "12F", "12G-H", "12I", "12J",
    >                 "12K-M", "12L", "12N", "12O", "12P", "12Q-R-T", "12S",
    >                 "24A-B-C-D", "24E-G", "24F", "24H", "24I-M", "24J", "24K",
    >                 "24L-N" ];
    >        one:= elm(2)^0;
    >        rank:= RankMat( elm(2) + one );
    >        pos:= Position( cand, [ order, rank ] );
    >        if not '-' in nams[ pos ] then
    >          return nams[ pos ];
    >        elif order = 12 then
    >          if nams[ pos ] = "12A-C-D" then
    >            trace:= Int( TraceMat( elm(5) ) );
    >            if   trace = 3 then return "12A";
    >            elif trace = 4 then return "12C";
    >            elif trace = 1 then return "12D";
    >            else Error( "wrong trace" );
    >            fi;
    >          else
    >            trace:= Int( TraceMat( elm(3) ) );
    >            if trace = 0 then
    >              if   nams[ pos ] = "12G-H" then return "12H";
    >              elif nams[ pos ] = "12K-M" then return "12K";
    >              else return "12Q";  # 12Q-R-T
    >              fi;
    >            elif trace = 1 then
    >              if   nams[ pos ] = "12G-H" then return "12G";
    >              elif nams[ pos ] = "12K-M" then return "12M";
    >              else return "12R";  # 12Q-R-T
    >              fi;
    >            elif nams[ pos ] = "12Q-R-T" then
    >              return "12T";
    >            else
    >              Error( "wrong trace" );
    >            fi;
    >          fi;
    >        elif order = 24 then
    >          if nams[ pos ] = "24I-M" then
    >            rank:= RankMat( elm(2)^2 + one );
    >            if   rank = 3986 then return "24I";
    >            elif rank = 3982 then return "24M";
    >            else Error( "wrong rank" );
    >            fi;
    >          elif nams[ pos ] = "24E-G" then
    >            rank:= RankMat( elm(2)^3 + one );
    >            if   rank = 3774 then return "24E";
    >            elif rank = 3778 then return "24G";
    >            else Error( "wrong rank" );
    >            fi;
    >          elif nams[ pos ] = "24L-N" then
    >            trace:= Int( TraceMat( elm(3) ) );
    >            if   trace = 1 then return "24L";
    >            elif trace = 2 then return "24N";
    >            else Error( "wrong trace" );
    >            fi;
    >          else   # 24A-B-C-D"
    >            trace:= Int( TraceMat( elm(5) ) );
    >            if   trace = 3 then return "24A";
    >            elif trace = 0 then return "24B";
    >            elif trace = 2 then return "24C";
    >            elif trace = 1 then return "24D";
    >            else Error( "wrong trace" );
    >            fi;
    >          fi;
    >        fi;
    >      elif order in [ 16, 20 ] then
    >        cand:= [ [ 16, 4072 ], [ 16, 4074 ], [ 16, 4094 ],
    >                 [ 20, 4114 ], [ 20, 4128 ], [ 20, 4132 ], [ 20, 4148 ],
    >                 [ 20, 4144 ], [ 20, 4138 ], [ 20, 4150 ] ];
    >        nams:= [ "16A-B", "16C-D-E-F", "16G-H",
    >                 "20A-B-C-D", "20E", "20F", "20G", "20H", "20I", "20J" ];
    >        one:= elm(2)^0;
    >        rank:= RankMat( elm(2) + one );
    >        pos:= Position( cand, [ order, rank ] );
    >        if not '-' in nams[ pos ] then
    >          return nams[ pos ];
    >        elif order = 20 then
    >          rank:= RankMat( elm(2)^2 + one );
    >          if   rank = 3908 then return "20B";
    >          elif rank <> 3896 then Error( "wrong rank" );
    >          else
    >            trace:= Int( TraceMat( elm(3) ) );
    >            if   trace = 2 then return "20A";
    >            elif trace = 0 then return "20C";
    >            elif trace = 1 then return "20D";
    >            else Error( "wrong trace" );
    >            fi;
    >          fi;
    >        else   # order = 16
    >          if nams[ pos ] = "16A-B" then
    >            trace:= Int( TraceMat( elm(3) ) );
    >            if   trace = 0 then return "16A";
    >            elif trace = 1 then return "16B";
    >            else Error( "wrong trace" );
    >            fi;
    >          elif nams[ pos ] = "16G-H" then
    >            trace:= Int( TraceMat( elm(3)^2 ) );
    >            if   trace = 1 then return "16G";
    >            elif trace = 2 then return "16H";
    >            else Error( "wrong trace" );
    >            fi;
    >          else   # 16C-D-E-F
    >            trace:= Int( TraceMat( elm(3) ) );
    >            if trace = 0 then   # We cannot distinguish 16D and 16F.
    >              return "16DF";
    >            elif trace = 2 then   # 16C-E
    >              one:= elm(2)^0;
    >              rank:= RankMat( elm(2)^2 + one );
    >              if   rank = 3780 then return "16C";
    >              elif rank = 3778 then return "16E";
    >              else Error( "wrong rank" );
    >              fi;
    >            else
    >              Error( "wrong trace" );
    >            fi;
    >          fi;
    >        fi;
    >      else Error( "wrong element order" );
    >      fi;
    >  end;;
\end{verbatim}

Elements of $\B$ to which the labels belong can be generated as follows.
The straight line program \texttt{"BG1-cycW1"} from~\cite{AGRv3}
computes generators of the maximally cyclic subgroups of $\B$.

\begin{verbatim}
    gap> cycprg:= AtlasProgram( "B", "cyclic" );;
    gap> cycprg.identifier;
    [ "B", "BG1-cycW1", 1 ]
    gap> cycprg.outputs;
    [ "12A", "12H", "12I", "12L", "12P", "12S", "12T", "16E", "16F", "16G", 
      "16H", "18A", "18B", "18D", "18F", "20B", "20C", "20H", "20I", "20J", 
      "24A", "24B", "24C", "24D", "24E", "24F", "24H", "24I", "24J", "24K", 
      "24L", "24M", "24N", "25A", "26B", "27A", "28A", "28C", "28D", "28E", 
      "30A", "30B", "30C", "30E", "30G-H", "31A-B", "32A-B", "32C-D", "34A", 
      "34BC", "36A", "36B", "36C", "38A", "39A", "40A", "40B", "40C", "40D", 
      "40E", "42A", "42B", "42C", "44A", "46AB", "47AB", "48A", "48B", "52A", 
      "55A", "56AB", "60A", "60B", "60C", "66A", "70A" ]
\end{verbatim}

The remaining representatives are obtained as suitable powers of them.
The following list encodes the definition of these powers,
see~\cite[Table~2]{BMverify}.

\begin{verbatim}
    gap> DefinitionsViaPowerMaps:= [
    >     [ "70A", 2, "35A" ], [ "66A", 2, "33A" ], [ "60A", 2, "30D" ],
    >     [ "60C", 2, "30F" ], [ "56AB", 2, "28B" ], [ "52A", 2, "26A" ],
    >     [ "48B", 2, "24G" ], [ "46AB", 2, "23AB" ], [ "44A", 2, "22B" ],
    >     [ "66A", 3, "22A" ], [ "42C", 2, "21A" ], [ "40E", 2, "20G" ],
    >     [ "40D", 2, "20F" ], [ "60B", 3, "20E" ], [ "60A", 3, "20A" ],
    >     [ "40C", 2, "20D" ], [ "38A", 2, "19A" ], [ "36C", 2, "18E" ],
    >     [ "36B", 2, "18C" ], [ "34A", 2, "17A" ], [ "32CD", 2, "16DF" ],
    >     [ "32AB", 2, "16C" ], [ "48B", 3, "16B" ], [ "48A", 3, "16A" ],
    >     [ "30F", 2, "15B" ], [ "30A", 2, "15A" ], [ "28E", 2, "14E" ],
    >     [ "28A", 2, "14D" ], [ "42A", 3, "14A" ], [ "42B", 3, "14B" ],
    >     [ "42C", 3, "14C" ], [ "26A", 2, "13A" ], [ "24N", 2, "12R" ],
    >     [ "24M", 2, "12O" ], [ "24L", 2, "12Q" ], [ "24K", 2, "12M" ],
    >     [ "24J", 2, "12J" ], [ "24H", 2, "12F" ], [ "24G", 2, "12G" ],
    >     [ "24D", 2, "12D" ], [ "36C", 3, "12N" ], [ "36B", 3, "12K" ],
    >     [ "36A", 3, "12B" ], [ "60A", 5, "12C" ], [ "60B", 5, "12E" ],
    >     [ "22B", 2, "11A" ], [ "20J", 2, "10F" ], [ "20I", 2, "10D" ],
    >     [ "20H", 2, "10C" ], [ "20F", 2, "10B" ], [ "30A", 3, "10A" ],
    >     [ "30E", 3, "10E" ], [ "18F", 2, "9B" ], [ "18E", 2, "9A" ],
    >     [ "16H", 2, "8M" ], [ "16G", 2, "8K" ], [ "16DF", 2, "8H" ],
    >     [ "16E", 2, "8D" ], [ "24J", 3, "8J" ], [ "24M", 3, "8I" ],
    >     [ "24I", 3, "8G" ], [ "24K", 3, "8F" ], [ "24C", 3, "8E" ],
    >     [ "24B", 3, "8C" ], [ "24A", 3, "8B" ], [ "24E", 3, "8A" ],
    >     [ "24N", 3, "8N" ], [ "40D", 5, "8L" ], [ "14D", 2, "7A" ],
    >     [ "12T", 2, "6K" ], [ "12S", 2, "6J" ], [ "12R", 2, "6I" ],
    >     [ "12P", 2, "6H" ], [ "12O", 2, "6G" ], [ "12I", 2, "6C" ],
    >     [ "18A", 3, "6D" ], [ "30B", 5, "6A" ], [ "30A", 5, "6B" ],
    >     [ "30E", 5, "6E" ], [ "30C", 5, "6F" ], [ "12C", 3, "4A" ],
    >     [ "10F", 2, "5B" ], [ "10B", 2, "5A" ], [ "8N", 2, "4J" ],
    >     [ "8M", 2, "4H" ], [ "8L", 2, "4G" ], [ "8J", 2, "4E" ],
    >     [ "8I", 2, "4F" ], [ "8H", 2, "4C" ], [ "8E", 2, "4B" ],
    >     [ "12E", 3, "4D" ], [ "12T", 3, "4I" ], [ "6K", 2, "3B" ],
    >     [ "6A", 2, "3A" ], [ "4J", 2, "2D" ], [ "4I", 2, "2C" ],
    >     [ "4A", 2, "2B" ], [ "6A", 3, "2A" ], [ "2B", 2, "1A" ],
    >    ];;
\end{verbatim}

The following function takes a label and the straight line program data
shown above, and returns a straight line program for computing an
element for the given label from standard generators of $\B$.

\begin{verbatim}
    gap> SLPForClassName:= function( nam, cycslp, outputnames )
    >      local pos, rule;
    > 
    >      pos:= Position( outputnames, nam );
    >      if pos <> fail then
    >        return RestrictOutputsOfSLP( cycslp.program, pos );
    >      fi;
    > 
    >      rule:= First( DefinitionsViaPowerMaps, x -> x[3] = nam );
    >      if rule = fail then
    >        Error( "'nam' is not an admiss. name for a cyclic subgroup of B" );
    >      fi;
    > 
    >      return CompositionOfStraightLinePrograms(
    >                 StraightLineProgram( [ [ 1, rule[2] ] ], 1 ),
    >                 SLPForClassName( rule[1], cycslp, outputnames ) );
    > end;;
\end{verbatim}

Let us verify that \texttt{IdentifyClassName} computes the claimed labels.

\begin{verbatim}
    gap> outputnames:= List( cycprg.outputs,
    >                        x -> ReplacedString( x, "-", "" ) );;
    gap> outputnames:= List( outputnames,
    >                        x -> ReplacedString( x, "16F", "16DF" ) );;
    gap> labels:= Union( outputnames,
    >                    List( DefinitionsViaPowerMaps, x -> x[3] ) );;
    gap> for l in labels do
    >      slp:= SLPForClassName( l, cycprg, outputnames );
    >      id:= IdentifyClassName( gens_2, gens_3, gens_5, slp, fail );
    >      if id <> l then
    >        Print( "#E  problem with identification: ", id, " vs. ", l, "\n" );
    >      fi;
    >    od;
\end{verbatim}


As we get no outputs, the identification is correct.

For later use, we collect power map information for the labels.
In order to simplify later loops,
we sort the labels w.~r.~t.~increasing element order.

\begin{verbatim}
    gap> SortParallel( List( labels, x -> Int( Filtered( x, IsDigitChar ) ) ),
    >                  labels );
    gap> powerinfo:= [];;
    gap> for l in labels do
    >      slp:= SLPForClassName( l, cycprg, outputnames );
    >      ord:= Int( Filtered( l, IsDigitChar ) );
    >      pow:= [];
    >      if not ( IsPrimeInt( ord ) or ord = 1 ) then
    >        for p in Set( Factors( ord ) ) do
    >          powerslp:= CompositionOfStraightLinePrograms(
    >                         StraightLineProgram( [ [ 1, p ] ], 1 ), slp );
    >          id:= IdentifyClassName( gens_2, gens_3, gens_5, powerslp,
    >                                  ord / p );
    >          Add( pow, [ p, id ] );
    >        od;
    >      fi;
    >      Add( powerinfo, [ l, pow ] );
    >    od;
\end{verbatim}

\section{Centralizers of elements of prime order}\label{centralizers_prime}

We document part of the computations needed in~\cite[Section~5]{BMverify}.

We know from \cite{Str76b} that $\B$ has exactly four classes of involutions,
whose normalizers in $\B$ have the following properties.

\begin{itemize}
\item
    The normalizer $H$ of a \texttt{2A} involution has the structure
    $2.{}^2E_6(2).2$, such that $H$ is an extension of its derived subgroup
    $H'$ by a field automorphism of order two, see~\cite[p.~505]{Str76b}.
    This implies that $H$ is a \emph{split} extension of $H'$,
    and this means that the character table of $H$ is the one that is shown
    in the {\ATLAS} of Finite Groups and that has the identifier
    \texttt{"2.2E6(2).2"} in {\GAP}'s library of character tables;
    note that the isoclinic variant of $H$ is a \emph{non-split} extension
    of $2.{}^2E_6(2)$.
\item
    The normalizer $C$ of a \texttt{2B} involution has the structure
    $2^{1+22}.Co_2$.
    Its character table has been constructed in Section~\ref{table_c2b}.
\item
    The normalizer of a \texttt{2C} involution is a subdirect product
    of $F_4(2).2$ and a dihedral group $D_8$ of order eight,
    see~\cite[La.~3.1, 5.6]{Str76b}.
    The character table of this group can easily be constructed
    character-theoretically from the known character tables of
    $F_4(2)$, $F_4(2).2$, $2^2$, and $D_8$.
\item
    The normalizer of a \texttt{2D} involution
    has order $11\,689\,182\,992\,793\,600$ and is contained in the
    normalizer in $\B$ of an elementary abelian group of order $2^8$.
    The character table of the latter normalizer has been computed from
    a subgroup of $\B$, see Appendix~\ref{2Dnormalizer}.
    This table has a unique class with centralizer order equal to the
    order of the normalizer of a \texttt{2D} element in $\B$.
\end{itemize}

Let $H \cong 2.{}^2E_6(2).2$ be a \texttt{2A} centralizer in $\B$.
The involution classes in $H$ are as stated
in the table in~\cite[Section~5.1]{BMverify}.

\begin{verbatim}
    gap> h:= CharacterTable( "2.2E6(2).2" );;
    gap> invpos:= Positions( OrdersClassRepresentatives( h ), 2 );
    [ 2, 3, 4, 5, 6, 7, 175, 176, 177, 178 ]
    gap> ClassNames( h ){ invpos };
    [ "2a", "2b", "2c", "2d", "2e", "2f", "2g", "2h", "2i", "2j" ]
    gap> AtlasClassNames( h ){ invpos };
    [ "1A_1", "2A_0", "2A_1", "2B_0", "2B_1", "2C_0", "2D_0", "2D_1", "2E_0", 
      "2E_1" ]
\end{verbatim}

(The subscripts $0$ and $1$ that appear above
are denoted by signs $+$ and $-$ in~\cite{BMverify}.)


Let us try to compute the necessary information about
the permutation action of $\B$ on the cosets of $H$,
restricted to $H$.
We calculate the first three transitive constituents
of the permutation character as listed in\cite[Section~5.2]{BMverify}.

The point stabilizer of the action of $H$ on the first orbit contains
the centre $Z(H)$ of $H$, thus we may perform the computations with
$H/Z(H)$.
If we assume that the rank of the permutation character is $5$ then
we can compute the possible degrees of the irreducible constituents
combinatorially, as follows.

\begin{verbatim}
    gap> f:= CharacterTable( "2E6(2).2" );;
    gap> index1:= 3968055;;
    gap> constit:= Filtered( Irr( f ), x -> x[1] <= index1 );;
    gap> degrees:= Set( constit, x -> x[1] );
    [ 1, 1938, 48620, 554268, 815100, 1322685, 1828332, 2089164, 2909907, 2956096 
     ]
    gap> lincom:= Filtered( UnorderedTuples( degrees, 5 ),
    >                       x -> Sum( x ) = index1 );
    [ [ 1, 1938, 48620, 1828332, 2089164 ] ]
\end{verbatim}

The degrees are uniquely determined.
Now we compute which sums of irreducibles of these degrees
have nonnegative values.

\begin{verbatim}
    gap> degrees:= lincom[1];
    [ 1, 1938, 48620, 1828332, 2089164 ]
    gap> constit:= List( degrees, d -> Filtered( constit, x -> x[1] = d ) );;
    gap> cand1:= List( Cartesian( constit ), Sum );;
    gap> cand1:= Filtered( cand1, x -> ForAll( x, y -> y >= 0 ) );;
    gap> List( ConstituentsOfCharacter( cand1[1] ),
    >          x -> Position( Irr( f ), x ) );
    [ 1, 3, 5, 13, 15 ]
\end{verbatim}

Thus this permutation character is uniquely determined.
Alternatively, we can ask {\GAP} to compute the possible permutation
characters of the given degree, and get the same result
(without a priori knowledge about the rank of the permutation action).

\begin{verbatim}
    gap> PermComb( f, rec( degree:= index1 ) ) = cand1;
    true
\end{verbatim}

We compute the permutation character of the action on the second orbit
in two steps.
First we induce the trivial character of $F_4(2)$ to $H$ and then
we compute the unique subcharacter of this character that has the right
degree and only nonnegative values.

\begin{verbatim}
    gap> u:= CharacterTable( "F4(2)" );;
    gap> ufush:= PossibleClassFusions( u, h );;
    gap> ind:= Set( ufush,
    >               map -> InducedClassFunctionsByFusionMap( u, h,
    >                          [ TrivialCharacter( u ) ], map )[1] );;
    gap> Length( ind );
    1
    gap> const:= ConstituentsOfCharacter( ind[1] );;
    gap> Sum( const ) = ind[1];
    true
    gap> sub:= List( Cartesian( List( const, x -> [ Zero( x ), x ] ) ), Sum );;
    gap> cand2:= Filtered( sub,
    >                      x -> x[1] = ind[1][1] / 4 and Minimum( x ) >= 0 );;
    gap> Length( cand2 );
    1
    gap> List( ConstituentsOfCharacter( cand2[1] ),
    >          x -> Position( Irr( h ), x ) );
    [ 1, 5, 17, 24 ]
\end{verbatim}

The character table of the point stabilizer $Fi_{22}.2$
of the action of $H$ on the third orbit is available.
We compute the corresponding permutation character by inducing the
trivial character of $Fi_{22}.2$ to $H$.
Note that the class fusion from $Fi_{22}.2$ to $H$ is unique up to
the group automorphism of $H$ that multiplies the elements outside
the derived subgroup of $H$ by the central involution in $H$;
we know that the class \texttt{2D} of the point stabilizer lies in a class
of $H$ that fuses into the class \texttt{2A} of $\B$,
thus the first of the two fusions is the right one.

\begin{verbatim}
    gap> s:= CharacterTable( "Fi22.2" );;
    gap> sfush:= PossibleClassFusions( s, h );;
    gap> Length( sfush );
    2
    gap> Positions( OrdersClassRepresentatives( s ), 2 );
    [ 2, 3, 4, 60, 61, 62 ]
    gap> List( sfush, x -> x[60] );
    [ 175, 176 ]
    gap> cand3:= InducedClassFunctionsByFusionMap( s, h,
    >                [ TrivialCharacter( s ) ], sfush[1] );;
    gap> SortedList( List( ConstituentsOfCharacter( cand3[1] ),
    >                      x -> Position( Irr( h ), x ) ) );
    [ 1, 3, 5, 13, 17, 28, 49, 76, 190, 192, 196, 202, 210, 217 ]
\end{verbatim}

Next we compute the value of the permutation character of the action
on the fourth orbit listed in\cite[Section~5.2]{BMverify} on the
class \texttt{3C} of $H$.
The point stabilizer is $H_5 \cong 2^{1+20}.U_4(3).2^2$,
thus the subgroup $H_5 Z(H) / Z(H)$ of $H / Z(H) \cong {}^2E_6(2).2$
lies in a maximal subgroup of type $(2^{1+20}:U_6(2)).2$,
see~\cite[p.~191]{CCN85}.

Because the extension of $2^{1+20}$ by $U_6(2)$ splits,
we know that $H_5 Z(H) / Z(H)$ has $U_4(3)$ type subgroups,
thus $H_5$ has subgroups of one of the types $U_4(3)$, $2.U_4(3)$.
Computing possible class fusions from both possibilities to $H$,
we get that the class of elements of order $3$ in $U_4(3)$ or $2.U_4(3)$
that belongs to the \texttt{3A} class of $U_4(3)$
lies in the class \texttt{3C} of $H$,
and the other classes of $3$-elements lie in the classes \texttt{3A} or
\texttt{3B} of $H$, as claimed.

\begin{verbatim}
    gap> h:= CharacterTable( "2.2E6(2).2" );;
    gap> u:= CharacterTable( "U4(3)" );;
    gap> Positions( OrdersClassRepresentatives( h ), 3 );
    [ 8, 10, 12 ]
    gap> ufush:= PossibleClassFusions( u, h );;
    gap> 3pos:= Positions( OrdersClassRepresentatives( u ), 3 );
    [ 3, 4, 5, 6 ]
    gap> Set( ufush, x -> x{ 3pos } );
    [ [ 12, 8, 10, 10 ], [ 12, 10, 8, 10 ] ]
    gap> u:= CharacterTable( "2.U4(3)" );;
    gap> ufush:= PossibleClassFusions( u, h );;
    gap> 3pos:= Positions( OrdersClassRepresentatives( u ), 3 );
    [ 5, 7, 9, 11 ]
    gap> Set( ufush, x -> x{ 3pos } );
    [ [ 12, 8, 10, 10 ], [ 12, 10, 8, 10 ] ]
\end{verbatim}

The \texttt{3A} elements in $U_4(3)$ have centralizer order $2^3 \cdot 3^6$
in this group.
The centralizer order in $2^{20}.U_4(3)$ is $2^5 \cdot 3^6$,
since the fixed space of a \texttt{3A} element on the unique $20$ dimensional
irreducible module in characteristic $2$ has dimension $2$
--this can be read off from the fact that the Brauer character value on
the class \texttt{3A} is $-7$.

\begin{verbatim}
    gap> u:= CharacterTable( "U4(3)" ) mod 2;;
    gap> phi:= Filtered( Irr( u ), x -> x[1] = 20 );;
    gap> Display( u, rec( chars:= phi, powermap:= false ) );
    U4(3)mod2
    
         2  7  3  2  2  .  .  .  .  .  .  .  .
         3  6  6  5  5  4  .  .  .  3  3  3  3
         5  1  .  .  .  .  1  .  .  .  .  .  .
         7  1  .  .  .  .  .  1  1  .  .  .  .
    
           1a 3a 3b 3c 3d 5a 7a 7b 9a 9b 9c 9d
    
    Y.1    20 -7  2  2  2  . -1 -1 -1 -1 -1 -1
\end{verbatim}

Now the centralizer order gets doubled in the central extension to
$2^{1+20}.U_4(3)$,
and the two upwards extensions cannot fuse the \texttt{3A} class with
another class, thus the centralizer order is again doubled in each case,
which means that $|C_{H_5}(\texttt{3A})| = 2^8 \cdot 3^6$.

The permutation character of the action of $\B$ on the cosets of $H$
has the value $1\,620$ on the class \texttt{3C} of $H$.

\begin{verbatim}
    gap> 3pos:= Positions( OrdersClassRepresentatives( f ), 3 );;
    gap> val3C:= 1 + cand1[1][ 3pos[3] ];;
    gap> 3pos:= Positions( OrdersClassRepresentatives( h ), 3 );;
    gap> val3C:= val3C + cand2[1][ 3pos[3] ] + cand3[1][ 3pos[3] ];;
    gap> val3C:= val3C + SizesCentralizers( h )[ 3pos[3] ] / ( 2^8 * 3^6 );
    1620
\end{verbatim}

Thus we have computed $|C_{\B}(\texttt{3B})| = 2^{13} \cdot 3^{13} \cdot 5$.

\begin{verbatim}
    gap> Collected( Factors( val3C * SizesCentralizers( h )[ 3pos[3] ] ) );
    [ [ 2, 13 ], [ 3, 13 ], [ 5, 1 ] ]
\end{verbatim}

\section{The character table of $2^{1+22}.Co_2$}\label{table_c2b}

Let $z$ be an involution in $\B$ whose class is called \texttt{2B} in the
{\ATLAS} of Finite Groups.
The centralizer $C$ of $z$ in $\B$ has the structure $2^{1+22}.Co_2$,
the construction of its character table is described in~\cite{Pah07},
and this table is available in {\GAP}.
However, that paper assumes the knowledge of the character table of $\B$,
hence we are not allowed to use the known character table of $C$
in the verification of the character table of $\B$.

In this section, we recompute the character table of $C$, as follows.
We start with the three certified matrix representations of $\B$
in characteristic $2$, $3$, and $5$,
and with the straight line program for restricting these representations
to a \texttt{2B} centralizer.

First we standardize the generators of the subgroup such that
the known straight line program for computing class representatives of $Co_2$
can be applied.
Next we compute a permutation representation of degree $4600$ for the
factor group $C / \langle z \rangle$.
Using this representation,
we construct a straight line program that computes class representatives
of $C / \langle Z \rangle$ from the images of the given generators of $C$
under the natural epimorphism.
Applying this straight line program to the restrictions of the
given representations of $\B$ and then computing the class labels in $\B$
yields a ``preliminary class fusion´´ from $C$ to $\B$.
Furthermore,
the given matrix representations of $C$ in characteristic $3$ and $5$
have a unique faithful constituent,
which lifts to the unique ordinary irreducible character of degree $2048$
of $C$.
Hence the Brauer characters of the two representations yield most of the
values of this ordinary character.
Together with the character table of $C / \langle Z \rangle$ that can be
computed by {\MAGMA}~\cite{Magma} from the permutation representation,
this information suffices to complete the character table of $C$.

Using the given degree $4371$ matrix representation of $\B$
over the field with three elements
(verified as described in Section~\ref{pres_B})
and the description from~\cite{AGRv3} how to restrict this representation
to (a conjugate of) $C$,
we compute generators of $C$.

\begin{verbatim}
    gap> slp_maxes_2:= AtlasProgram( "B", "maxes", 2 );;
    gap> cgens_3:= ResultOfStraightLineProgram( slp_maxes_2.program,
    >                  GeneratorsOfGroup( b_3 ) );;
\end{verbatim}

The composition factors of the $4371$ dimensional module for $C$ have the
dimensions $23$, $2300$, and $2048$, respectively.
The kernels of the actions of $C$ on these factors will turn out to have
the orders $2^{23}$, $2$, and $1$, respectively.

\begin{verbatim}
    gap> m:= GModuleByMats( cgens_3, GF(3) );;
    gap> cf_3:= MTX.CompositionFactors( m );;
    gap> SortParallel( List( cf_3, x -> x.dimension ), cf_3 );
    gap> List( cf_3, x -> x.dimension );
    [ 23, 2048, 2300 ]
\end{verbatim}

We use the action on the $23$-dimensional module
to find words in terms of the generators of $C$
that act on this module as \emph{standard generators} of $Co_2$,
as defined in~\cite{Wil96}.

For that, we find two elements $a$, $b$ that generate $Co_2$ and
lie in the conjugacy classes \texttt{2A} and \texttt{5A}, respectively,
such that the product $a b$ has order $28$.

%
%

Let us call the generators in the $23$-dimensional composition factor
$x$ and $y$, and set $c = y (y x)^3$.
Then the elements $y^{12}$, $c^{-1} ((y^4 x)^4) c$ are standard generators
of $Co_2$, see Appendix~\ref{slp_co2}.
The following straight line program computes these elements
when it is applied to $x$ and $y$.

\begin{verbatim}
    gap> slp_co2:= StraightLineProgram( [
    >     [ 2, 1, 1, 1 ], [ 2, 2 ], [ 4, 1, 1, 1 ], [ 4, 2 ],
    >     [ 6, 1, 1, 1 ], [ 7, 4 ], [ 3, 2 ], [ 5, 1, 9, 1 ],
    >     [ 10, -1 ], [ 6, 3 ], [ 11, 1, 8, 1 ], [ 13, 1, 10, 1 ],
    >     [ [ 12, 1 ], [ 14, 1 ] ] ], 2 );;
    gap> f:= FreeGroup( "x", "y" );;  x:= f.1;;  y:= f.2;;
    gap> words:= ResultOfStraightLineProgram( slp_co2, [ x, y ] );;
    gap> words = [ y^12, ((y^4*x)^4)^(y*(y*x)^3) ];
    true
    gap> co2gens:= cf_3[1].generators;;
    gap> co2stdgens:= ResultOfStraightLineProgram( slp_co2, co2gens );;
\end{verbatim}

Next we find an orbit of length $4600$ under the action of $C$
on the $2300$-dimensional module.
This will allow us to represent $C/\langle z \rangle$ as a permutation group
$P$, say, of degree $4600$.

Note that $Co_2$ contains a maximal subgroup of the structure $U_6(2).2$,
of index $2300$,
and that there is a $2^{21}.U_6(2).2$ type subgroup of $C$ that fixes
a $1$-dimensional subspace in the given $2300$-dimensional representation
of $C$ over the field with three elements.
We can compute generators for (a sufficiently large subgroup of)
$2^{22}.U_6(2).2$ by applying the straight line program for computing
generators of a $U_6(2).2$ type subgroup from standard generators of $Co_2$.
Then we find a vector in the $1$-dimensional subspace
as the common fixed vector of squares of commutators in this subgroup.


\begin{verbatim}
    gap> fgens:= cf_3[3].generators;;
    gap> fstdgens:= ResultOfStraightLineProgram( slp_co2, fgens );;
    gap> slp_co2m1:= AtlasProgram( "Co2", "maxes", 1 );;
    gap> ugens:= ResultOfStraightLineProgram( slp_co2m1.program, fstdgens );;
    gap> one:= ugens[1]^0;;
    gap> comm:= Comm( ugens[1], ugens[2] );;
    gap> Order( comm );
    12
    gap> pow:= comm^2;;
    gap> mats:= List( [ pow, pow^ugens[1] ], x -> x - one );;
    gap> nsp:= List( mats, NullspaceMat );;
    gap> List( nsp, Length );
    [ 434, 434 ]
    gap> si:= SumIntersectionMat( nsp[1], nsp[2] );;
    gap> Length( si[2] );
    1
    gap> orb:= Orbit( Group( fstdgens ), si[2][1] );;
    gap> Length( orb );
    4600
    gap> orb:= SortedList( orb );;
    gap> stdperms:= List( fstdgens, x -> Permutation( x, orb ) );;
    gap> List( stdperms, Order );
    [ 2, 5 ]
\end{verbatim}

In $P$, we compute a basis for the elementary abelian
normal subgroup $N$ of order $2^{22}$,
and the following words for the basis vectors
in terms of the standard generators of $Co_2$,
see Appendix~\ref{slp_kernel}.


\[
\begin{array}[t]{llll}
   (b^2 a)^9,   & (b a b)^9,     & (a b^2)^9,       & (b^2 a b)^9, \\
   (b a b^2)^9, & ((a b)^2 a)^9, & (a b (b a)^2)^9, & ((a b)^2 b a)^9, \\
   (b^3 (b a)^2)^4, & (b (b^2 a)^2)^4, & (b^2 a b^3 a)^4, & (b a b^4 a)^4, \\
   (b^2 (b a)^2 b)^4, & ((b^2 a)^2 b)^4, & (b a b^3 a b)^4, &
   (b (b a)^2 b^2)^4, \\
   ((b a b)^2 b)^4, & ((b^2 a)^2 b a)^{12}, & (b (b a)^2 b^2 a)^{12}, &
   ((b a b)^2 b a)^{12}, \\
   ((b a b)^2 a b)^{12}, & ((b^2 a)^2 b a b a)^{12} & &
\end{array}
\]

The following straight line program computes these basis vectors
when it is applied to $x$ and $y$.

\begin{verbatim}
    gap> slp_ker:= StraightLineProgram( [
    >     [ 1, 1, 2, 1 ], [ 2, 1, 1, 1 ], [ 2, 2 ], [ 4, 1, 2, 1 ],
    >     [ 2, 1, 4, 1 ], [ 3, 1, 2, 1 ], [ 3, 2 ], [ 4, 2 ], [ 6, 2 ],
    >     [ 7, 2 ], [ 2, 1, 10, 1 ], [ 2, 1, 6, 1 ], [ 10, 1, 2, 1 ],
    >     [ 5, 1, 3, 1 ], [ 4, 1, 5, 1 ], [ 9, 1, 1, 1 ], [ 3, 1, 10, 1 ],
    >     [ 9, 1, 4, 1 ], [ 5, 1, 13, 1 ], [ 2, 1, 12, 1 ], [ 14, 1, 7, 1 ],
    >     [ 17, 1, 7, 1 ], [ 5, 1, 15, 1 ], [ 12, 1, 2, 1 ], [ 6, 1, 14, 1 ],
    >     [ 13, 1, 5, 1 ], [ 11, 1, 2, 1 ], [ 12, 1, 4, 1 ], [ 13, 1, 7, 1 ],
    >     [ 11, 1, 4, 1 ], [ 11, 1, 3, 1 ], [ 12, 1, 10, 1 ],
    >     [ [ 7, 9 ], [ 6, 9 ], [ 8, 9 ], [ 16, 9 ], [ 17, 9 ], [ 18, 9 ],
    >       [ 19, 9 ], [ 20, 9 ], [ 21, 4 ], [ 22, 4 ], [ 23, 4 ], [ 24, 4 ],
    >       [ 25, 4 ], [ 26, 4 ], [ 27, 4 ], [ 28, 4 ], [ 29, 4 ], [ 30, 12 ],
    >       [ 31, 12 ], [ 32, 12 ], [ 33, 12 ], [ 34, 12 ] ] ], 2 );;
    gap> f:= FreeGroup( "a", "b" );;  a:= f.1;;  b:= f.2;;
    gap> ResultOfStraightLineProgram( slp_ker, [ a, b ] );
    [ (b^2*a)^9, (b*a*b)^9, (a*b^2)^9, (b^2*a*b)^9, (b*a*b^2)^9, ((a*b)^2*a)^9, 
      (a*b*(b*a)^2)^9, ((a*b)^2*b*a)^9, (b^3*(b*a)^2)^4, (b*(b^2*a)^2)^4, 
      (b^2*a*b^3*a)^4, (b*a*b^4*a)^4, (b^2*(b*a)^2*b)^4, ((b^2*a)^2*b)^4, 
      (b*a*b^3*a*b)^4, (b*(b*a)^2*b^2)^4, ((b*a*b)^2*b)^4, ((b^2*a)^2*b*a)^12, 
      (b*(b*a)^2*b^2*a)^12, ((b*a*b)^2*b*a)^12, ((b*a*b)^2*a*b)^12, 
      ((b^2*a)^2*b*a*b*a)^12 ]
\end{verbatim}

the straight line program that is available in~\cite{AGRv3}
computes class representatives of $Co_2$ from standard generators.

\begin{verbatim}
    gap> slp_co2classreps:= AtlasProgram( "Co2", "classes" );;
\end{verbatim}

Hence we can describe representatives of the $388$ classes of $P$
as products of the class representatives of $Co_2$
and suitable elements of $N$.
The necessary computations are described in Appendix~\ref{slp_classreps}.

The list \texttt{classrepsinfo} contains $60$ entries;
the $i$-th entry describes the preimage classes of the $i$-th class
of $Co_2$,
by listing the positions of those basis vectors that can be multiplied
with the $i$-th output of \texttt{slp\_co2classreps} in order to get
the class representatives in question.

(The ordering of the representatives fits to the ordering of the classes
in the relevant factor of the library character table of $C$.
The entry \texttt{[ 0 ]} means that the identity matrix is taken instead
of the representative.)

\begin{verbatim}
    gap> classrepsinfo:= [
    >  [  "1A", [ [ 0 ], [ 3, 4, 22 ], [  ], [ 8 ], [ 2 ], [ 1 ] ] ], 
    >  [  "2A", [ [ 1, 6, 8, 9, 11 ], [ 12 ], [ 1, 4, 9, 10, 19 ],
    >             [ 5, 7 ], [ 1, 9 ], [ 1, 16 ], [  ], [ 1 ], [ 5 ],
    >             [ 1, 5, 18 ], [ 3 ] ] ], 
    >  [  "2B", [ [ 1, 3, 4, 17, 20 ], [ 2, 14, 17 ], [ 2, 10, 22 ],
    >             [ 14, 21 ], [  ], [ 9 ], [ 1, 5 ], [ 10, 11 ], [ 6, 12 ],
    >             [ 1 ], [ 2 ] ] ], 
    >  [  "2C", [ [ 4, 11, 13, 20 ], [ 1, 10, 12, 17 ], [ 3, 18, 21 ],
    >             [ 8, 13, 17 ], [ 16, 22 ], [ 4, 5, 15 ], [ 8, 13 ],
    >             [ 1, 10 ], [ 1, 21 ], [ 3, 9, 18, 21 ], [ 8 ],
    >             [ 1, 10, 12 ], [ 4 ], [ 21 ], [ 6 ], [ 1 ], [  ] ] ],
    >  [  "3A", [ [  ], [ 21 ], [ 1 ] ] ], 
    >  [  "3B", [ [ 17, 20 ], [ 17 ], [ 1, 5, 17 ], [ 3, 18 ], [ 1 ], [ 10 ],
    >             [ 12 ], [ 3 ], [  ] ] ],
    >  [  "4A", [ [ 11 ], [  ], [ 1 ], [ 1, 12 ], [ 2 ] ] ], 
    >  [  "4B", [ [ 2, 7, 19, 20 ], [ 2, 12 ], [ 9, 21 ], [ 5 ], [  ],
    >             [ 14 ], [ 9, 20 ], [ 2, 14 ], [ 6 ], [ 1 ], [ 8 ], [ 2 ] ] ], 
    >  [  "4C", [ [ 3, 7, 11 ], [ 18, 22 ], [ 4, 6 ], [ 21 ], [ 2, 17, 18 ],
    >             [ 1, 14 ], [ 5, 12, 17 ], [ 14 ], [ 3, 4 ], [ 3 ], [ 10 ],
    >             [ 2 ], [ 1 ], [  ] ] ], 
    >  [  "4D", [ [ 20 ], [ 9, 22 ], [ 8, 16 ], [ 6 ], [ 6, 17 ], [  ],
    >             [ 1 ] ] ], 
    >  [  "4E", [ [ 11, 19 ], [ 3, 7 ], [ 1, 22 ], [ 1, 2, 22 ], [ 1, 10 ],
    >             [ 3, 18 ], [  ], [ 5, 20 ], [ 2 ], [ 1, 15 ], [ 1 ], [ 3 ],
    >             [ 8 ] ] ], 
    >  [  "4F", [ [ 3, 4 ], [ 3, 19, 21 ], [ 4, 16 ], [ 4, 13 ], [ 3, 7 ],
    >             [ 4 ], [ 3 ], [ 4, 17 ], [ 13 ], [ 1, 7 ], [  ], [ 1, 9 ],
    >             [ 4, 14 ], [ 17 ], [ 19 ], [ 1 ], [ 1, 4 ], [ 18 ],
    >             [ 11 ] ] ], 
    >  [  "4G", [ [ 10 ], [ 1, 4 ], [ 8 ], [ 20 ], [  ], [ 16 ], [ 1 ],
    >             [ 2 ] ] ], 
    >  [  "5A", [ [  ], [ 1 ] ] ], 
    >  [  "5B", [ [  ], [ 1 ], [ 2 ], [ 19 ], [ 15, 19 ], [ 4 ], [ 8 ] ] ], 
    >  [  "6A", [ [ 2 ], [  ], [ 1 ] ] ], 
    >  [  "6B", [ [  ], [ 21 ], [ 1, 11 ], [ 1 ], [ 10 ] ] ], 
    >  [  "6C", [ [ 8 ], [ 1, 14 ], [  ], [ 9 ], [ 1 ], [ 10 ], [ 2 ],
    >             [ 3 ] ] ], 
    >  [  "6D", [ [  ], [ 19 ], [ 12, 17 ], [ 1, 11 ], [ 7 ], [ 4, 10 ], [ 1 ], 
    >             [ 2, 3 ], [ 1, 7 ], [ 10 ], [ 3 ] ] ], 
    >  [  "6E", [ [ 9 ], [ 1, 8 ], [ 2 ], [ 22 ], [  ], [ 2, 3 ], [ 1 ],
    >             [ 3, 4 ], [ 3 ], [ 13 ], [ 7 ], [ 6 ] ] ], 
    >  [  "6F", [ [ 5, 10 ], [ 1, 2, 4 ], [ 4, 13 ], [ 10, 18 ], [ 4, 12 ],
    >             [ 4, 9 ], [ 21 ], [ 4 ], [ 8 ], [ 1, 10 ], [ 1, 8, 18 ],
    >             [ 10 ], [ 6 ], [ 1 ], [  ], [ 3 ] ] ], 
    >  [  "7A", [ [ 1, 2 ], [ 2 ], [  ], [ 7 ], [ 1, 10 ], [ 1 ], [ 11 ] ] ], 
    >  [  "8A", [ [ 2, 18 ], [  ], [ 2 ], [ 9 ], [ 7 ], [ 1 ] ] ], 
    >  [  "8B", [ [ 8, 17 ], [ 1, 7 ], [ 11, 21 ], [ 2, 8 ], [  ], [ 1 ],
    >             [ 2, 12 ], [ 7 ], [ 5 ], [ 2 ] ] ], 
    >  [  "8C", [ [ 1, 8 ], [ 13 ], [ 9 ], [  ], [ 1 ] ] ], 
    >  [  "8D", [ [ 1, 2 ], [ 7 ], [ 1 ], [ 6 ], [  ] ] ], 
    >  [  "8E", [ [ 2, 12 ], [ 3, 22 ], [ 4, 22 ], [ 4, 12 ], [ 1, 11 ],
    >             [ 2, 22 ], [ 6, 12 ], [ 9 ], [ 2, 10 ], [ 20 ], [ 1, 9 ],
    >             [ 11 ], [ 10 ], [ 1 ], [  ], [ 7 ] ] ], 
    >  [  "8F", [ [ 1, 10 ], [ 18 ], [ 4 ], [ 9 ], [ 20 ], [  ], [ 6 ],
    >             [ 1 ] ] ], 
    >  [  "9A", [ [  ], [ 4 ], [ 3 ], [ 5 ] ] ],
    >  [ "10A", [ [  ], [ 1 ] ] ], 
    >  [ "10B", [ [ 2 ], [ 12 ], [  ], [ 1 ], [ 9 ], [ 10 ], [ 14 ], [ 3 ] ] ], 
    >  [ "10C", [ [ 1, 20 ], [ 5 ], [ 1, 8 ], [ 12 ], [ 2 ], [  ], [ 10 ],
    >             [ 3 ], [ 1 ] ] ],
    >  [ "11A", [ [ 2 ], [ 1 ], [ 3 ], [  ] ] ], 
    >  [ "12A", [ [ 10 ], [  ], [ 1 ] ] ], 
    >  [ "12B", [ [  ], [ 14 ], [ 1 ], [ 3 ], [ 4 ], [ 10 ] ] ], 
    >  [ "12C", [ [ 19 ], [  ], [ 7 ], [ 1, 7 ], [ 1 ] ] ], 
    >  [ "12D", [ [ 1 ], [ 11 ], [ 8 ], [ 15 ], [  ], [ 2 ] ] ], 
    >  [ "12E", [ [ 3 ], [ 5 ], [  ] ] ], 
    >  [ "12F", [ [ 9 ], [ 1, 15 ], [ 5 ], [  ], [ 2 ], [ 12 ], [ 8 ],
    >             [ 1 ] ] ], 
    >  [ "12G", [ [  ], [ 1 ], [ 3 ] ] ], 
    >  [ "12H", [ [ 1, 4 ], [ 4, 14 ], [ 4, 9 ], [ 4 ], [ 14 ], [ 9 ], [ 1 ],
    >             [  ], [ 11 ], [ 18 ] ] ], 
    >  [ "14A", [ [ 1, 2 ], [ 1, 7 ], [ 10 ], [ 1 ], [ 19 ], [  ], [ 16 ] ] ], 
    >  [ "14B", [ [ 7 ], [  ], [ 1 ], [ 11 ] ] ], 
    >  [ "14C", [ [ 7 ], [  ], [ 1 ], [ 11 ] ] ], 
    >  [ "15A", [ [  ], [ 2 ], [ 3 ], [ 1 ] ] ],
    >  [ "15B", [ [  ] ] ], 
    >  [ "15C", [ [  ] ] ],
    >  [ "16A", [ [ 6 ], [ 11 ], [  ], [ 1 ] ] ], 
    >  [ "16B", [ [ 3 ], [  ], [ 6 ], [ 1 ] ] ], 
    >  [ "18A", [ [ 4 ], [  ], [ 5 ], [ 3 ] ] ], 
    >  [ "20A", [ [ 6 ], [ 1 ], [  ], [ 7 ] ] ], 
    >  [ "20B", [ [ 3 ], [ 1 ], [ 2 ], [  ] ] ],
    >  [ "23A", [ [  ] ] ], 
    >  [ "23B", [ [  ] ] ],
    >  [ "24A", [ [ 18 ], [ 7 ], [  ], [ 1 ] ] ], 
    >  [ "24B", [ [ 10 ], [  ], [ 1 ], [ 4 ] ] ], 
    >  [ "28A", [ [ 5 ], [  ], [ 10 ], [ 1 ] ] ], 
    >  [ "30A", [ [ 2 ], [  ], [ 1 ], [ 3 ] ] ],
    >  [ "30B", [ [  ] ] ], 
    >  [ "30C", [ [  ] ] ] ];;
\end{verbatim}

The {\GAP} code for turning this information into a straight line program
is shorter than the lines of this program,
hence we show this program.

\begin{verbatim}
    gap> create_classreps_slp:= function( classreps )
    >     local words, l, len, lines, kerneloffset, inputoffset, cache, k,
    >           found, pair, pos, diff, pos2, first, n, outputs, i, list;
    > 
    >     # Find words for the products of kernel generators that occur.
    >     words:= [];
    >     for l in Set( Concatenation( List( classreps, x -> x[2] ) ) ) do
    >       len:= Length( l );
    >       if 2 <= len then
    >         if not IsBound( words[ len ] ) then
    >           words[ len ]:= [];
    >         fi;
    >         Add( words[ len ], l );
    >       fi;
    >     od;
    > 
    >     lines:= [];
    >     kerneloffset:= 60;
    >     inputoffset:= 82;
    > 
    >     # We have to form all products of length 2 of kernel generators.
    >     cache:= [ [], [] ];
    >     for l in words[2] do
    >       Add( lines,
    >            [ l[1] + kerneloffset, 1, l[2] + kerneloffset, 1 ] );
    >       Add( cache[1], l );
    >       Add( cache[2], Length( lines ) + inputoffset );
    >     od;
    > 
    >     # For products of length at least 3, we may use known products
    >     # of length 2.  Longer matches are not considered.
    >     for k in [ 3 .. Length( words ) ] do
    >       for l in words[k] do
    >         found:= false;
    >         for pair in Combinations( l, 2 ) do
    >           pos:= Position( cache[1], pair );
    >           if pos <> fail then
    >             diff:= Difference( l, pair );
    >             if Length( diff ) = 1 then
    >               Add( lines,
    >                 [ cache[2][ pos ], 1, diff[1] + kerneloffset, 1 ] );
    >             else
    >               pos2:= Position( cache[1], diff );
    >               if pos2 <> fail then
    >                 Add( lines,
    >                   [ cache[2][ pos ], 1, cache[2][ pos2 ], 1 ] );
    >               else
    >                 first:= cache[2][ pos ];
    >                 for n in diff do
    >                   Add( lines, [ first, 1, n + kerneloffset, 1 ] );
    >                   first:= Length( lines ) + inputoffset;
    >                 od;
    >               fi;
    >             fi;
    >             Add( cache[1], l );
    >             Add( cache[2], Length( lines ) + inputoffset );
    >             found:= true;
    >             break;
    >           fi;
    >         od;
    >         if not found then
    >           first:= l[1] + kerneloffset;
    >           for n in l{ [ 2 .. Length( l ) ] } do
    >             Add( lines, [ first, 1, n + kerneloffset, 1 ] );
    >             first:= Length( lines ) + inputoffset;
    >           od;
    >           Add( cache[1], l );
    >           Add( cache[2], Length( lines ) + inputoffset );
    >         fi;
    >       od;
    >     od;
    > 
    >     outputs:= [];
    > 
    >     for i in [ 1 .. Length( classreps ) ] do
    >       list:= classreps[i][2];
    >       for l in list do
    >         if l = [ 0 ] then
    >           Add( outputs, [ kerneloffset + 1, 2 ] );
    >         elif l = [] then
    >           Add( outputs, [ i, 1 ] );
    >         elif Length( l ) = 1 then
    >           Add( lines, [ i, 1, l[1] + kerneloffset, 1 ] );
    >           Add( outputs, [ Length( lines ) + inputoffset, 1 ] );
    >         else
    >           # The words are already cached.
    >           pos:= Position( cache[1], l );
    >           Add( lines, [ i, 1, cache[2][ pos ], 1 ] );
    >           Add( outputs, [ Length( lines ) + inputoffset, 1 ] );
    >         fi;
    >       od;
    >     od;
    > 
    >     Add( lines, outputs );
    > 
    >     return StraightLineProgram( lines, inputoffset );
    > end;;
    gap> slp_classreps:= create_classreps_slp( classrepsinfo );;
\end{verbatim}

We compute the class representatives of $P$.

\begin{verbatim}
    gap> kerperms:= ResultOfStraightLineProgram( slp_ker, stdperms );;
    gap> co2classreps:= ResultOfStraightLineProgram(
    >        slp_co2classreps.program, stdperms );;
    gap> classreps:= ResultOfStraightLineProgram( slp_classreps,
    >        Concatenation( co2classreps, kerperms ) );;
\end{verbatim}

Now the {\MAGMA}~\cite{Magma} system is invoked
for computing the irreducible characters of $P$.
The function \texttt{CharacterTableComputedByMagma} guarantees that
the columns are indexed by the class representatives we have chosen.

\begin{verbatim}
    gap> g:= Group( stdperms );;
    gap> SetConjugacyClasses( g,
    >        List( classreps, x -> ConjugacyClass( g, x ) ) );
    gap> libcb2b:= CharacterTable( "BM2" );
    CharacterTable( "2^(1+22).Co2" )
    gap> cen:= ClassPositionsOfCentre( libcb2b );;
    gap> if CTblLib.IsMagmaAvailable() then
    >      mgmt:= CharacterTableComputedByMagma( g, "2^22.Co2-Magma" );
    >    else
    >      mgmt:= libcb2b / cen;  # this is a hack ...
    >    fi;
\end{verbatim}



Our next goal is to compute the character table of $C$,
together with the information about the correspondence of the
conjugacy classes in $C$ and $\B$.

In order to write down class representatives of $C$,
we evaluate the words for class representatives of
$C / \langle z \rangle$ in the generators of the faithful
$2048$-dimensional module for $C$, in characteristic $3$ and $5$.
For elements of order not divisible by $15$, we can compute the Brauer
character value in at least one of the two representations,
and interpret it as the value of the unique faithful irreducible ordinary
character of degree $2048$ of $C$.
Whenever this character value is nonzero, we know that the corresponding
class of $C / \langle z \rangle$ splits into two classes of $C$,
on which the values of this character differ by sign.
For classes where the character value is zero, it will turn out later
that no splitting occurs; for the moment, we leave this question open.

We get words in terms of the generators $a$ and $b$ of $C$ for the elements
in question, that is, for the class representatives of $C / \langle z \rangle$
and for products of some of them with the central involution of $C$.

We evaluate these words in the three $4371$-dimensional representations,
and run the identification program in order to assign the label of one of
the {\lq\lq preliminary conjugacy classes\rq\rq} of $\B$ to them.

Let us collect the necessary data, that is,
the class representatives of $C / \langle z \rangle$ in the restrictions
of the three representations of $\B$ to $C$
and in the two $2048$-dimensional representations of $C$,
in characteristics $3$ and $5$, respectively.

\begin{verbatim}
    gap> cgens_2:= ResultOfStraightLineProgram( slp_maxes_2.program,
    >                  GeneratorsOfGroup( b_2 ) );;
    gap> cgens_5:= ResultOfStraightLineProgram( slp_maxes_2.program,
    >                  GeneratorsOfGroup( b_5 ) );;
    gap> cgens_2_std:= ResultOfStraightLineProgram( slp_co2, cgens_2 );;
    gap> cgens_3_std:= ResultOfStraightLineProgram( slp_co2, cgens_3 );;
    gap> cgens_5_std:= ResultOfStraightLineProgram( slp_co2, cgens_5 );;
    gap> m:= GModuleByMats( cgens_5_std, GF(5) );;
    gap> cf_5:= MTX.CompositionFactors( m );;
    gap> SortParallel( List( cf_5, x -> x.dimension ), cf_5 );
    gap> List( cf_5, x -> x.dimension );
    [ 23, 2048, 2300 ]
    gap> inputsB:= List( [ cgens_2_std, cgens_3_std, cgens_5_std ],
    >      l -> Concatenation(
    >             ResultOfStraightLineProgram( slp_co2classreps.program, l ),
    >             ResultOfStraightLineProgram( slp_ker, l ) ) );;
    gap> cf3std:= ResultOfStraightLineProgram( slp_co2, cf_3[2].generators );;
    gap> inputs2048:= List( [ cf3std, cf_5[2].generators ],
    >      l -> Concatenation(
    >             ResultOfStraightLineProgram( slp_co2classreps.program, l ),
    >             ResultOfStraightLineProgram( slp_ker, l ) ) );;
\end{verbatim}


Next we need the central involution of $C$ in the five representations.
We are lucky, the preimage of the identity of $C / \langle z \rangle$
under the epimorphism from $C$ that is computed by the straight line
program \texttt{slp\_classreps} has order two.

\begin{verbatim}
    gap> slp:= RestrictOutputsOfSLP( slp_classreps, 1 );;
    gap> centralinv_2048:= List( inputs2048,
    >        l -> ResultOfStraightLineProgram( slp, l ) );;
    gap> List( centralinv_2048, Order );
    [ 2, 2 ]
    gap> centralinv_B:= List( inputsB,
    >        l -> ResultOfStraightLineProgram( slp, l ) );;
\end{verbatim}

Now we run over the class representatives of $C / \langle z \rangle$,
and collect the data in a record with the following components.

\texttt{preimages}:
    one or two preimage classes, depending on whether the class in question
    need not split or must split,

\texttt{fusionlabels}:
    for each class of $C / \langle z \rangle$,
    one or two labels of class names in $\B$,

\texttt{projcharacter}:
    the Brauer character value of the first preimage of the class
    in one of the $2048$-dimensional representations, if possible;
    otherwise \texttt{fail}.

\begin{verbatim}
    gap> cent_table:= rec( preimages:= [],
    >                      fusionlabels:= [],
    >                      projcharacter:= [] );;
    gap> for i in [ 1 .. Length( classreps ) ] do
    >      # identify the representative
    >      slp:= RestrictOutputsOfSLP( slp_classreps, i );
    >      id:= IdentifyClassName( inputsB[1], inputsB[2], inputsB[3], slp, fail );
    >      order:= Int( Filtered( id, IsDigitChar ) );
    > 
    >      if order mod 3 = 0 then
    >        if order mod 5 = 0 then
    >          # We cannot compute the Brauer character value.
    >          value:= fail;
    >        else
    >          value:= BrauerCharacterValue(
    >                      ResultOfStraightLineProgram( slp, inputs2048[2] ) );
    >        fi;
    >      else
    >        value:= BrauerCharacterValue(
    >                    ResultOfStraightLineProgram( slp, inputs2048[1] ) );
    >      fi;
    > 
    >      if value = 0 then
    >        # Assume no splitting.
    >        Add( cent_table.preimages, [ i ] );
    >        Add( cent_table.fusionlabels, [ id ] );
    >        Add( cent_table.projcharacter, value );
    >      else
    >        # Identify the class of the other preimage.
    >        mats:= List( inputsB,
    >                     l -> ResultOfStraightLineProgram( slp, l ) );
    >        id2:= IdentifyClassName( mats[1] * centralinv_B[1],
    >                            mats[2] * centralinv_B[2],
    >                            mats[3] * centralinv_B[3],
    >                            fail, fail );
    >        if value = fail then
    >          # no Brauer character value known
    >          Add( cent_table.preimages, [ i ] );
    >          Add( cent_table.fusionlabels, [ id, id2 ] );
    >          Add( cent_table.projcharacter, value );
    >        else
    >          # two preimage classes, take the positive value first
    >          Add( cent_table.preimages, [ i, i ] );
    >          if value > 0 then
    >            Add( cent_table.fusionlabels, [ id, id2 ] );
    >            Add( cent_table.projcharacter, value );
    >          else
    >            Add( cent_table.fusionlabels, [ id2, id ] );
    >            Add( cent_table.projcharacter, -value );
    >          fi;
    >        fi;
    >      fi;
    >    od;
\end{verbatim}

Let us compute the missing values for the faithful irreducible character
of $C$.
We know that these values are integers,
and the character values at the $p$-th powers are known,
for $p \in \{ 3, 5 \}$.
The value at the $p$-th power of an element $g$, say,
determines the congruence class of the value at $g$ modulo $p$,
thus we know the congruence classes of the missing values modulo $15$.
For each of the classes of $C$ where the character value is not known yet,
the class length is at least $|C| / 120$,
and a character value of absolute value $7$ or larger
on any of these classes would lead to a contribution of at least $7^2 / 120$
to the norm of the character.
However, the known character values contribute already more than
$1 - 7^2 / 120$ to the norm, hence the missing character values
are uniquely determined by their congruence class modulo $15$.

\begin{verbatim}
    gap> chi:= cent_table.projcharacter;;
    gap> failpos:= Positions( chi, fail );;
    gap> OrdersClassRepresentatives( mgmt ){ failpos };
    [ 15, 30, 30, 30, 15, 15, 30, 30, 30, 30, 30, 30 ]
    gap> SizesCentralizers( mgmt ){ failpos };
    [ 120, 120, 120, 120, 30, 30, 120, 120, 120, 120, 30, 30 ]
    gap> used:= 0;;
    gap> for i in[ 1 .. 388 ] do
    >      if chi[i] <> 0 and chi[i] <> fail then
    >        used:= used + 2 * SizesConjugacyClasses( mgmt )[i] * chi[i]^2;
    >      fi;
    >    od;
    gap> used / ( 2 * Size( mgmt ) ) > 1 - 7^2 / 120;
    true
\end{verbatim}

Thus we may complete the character, as follows.

\begin{verbatim}
    gap> for i in failpos do
    >      v:= ChineseRem( [ 3, 5 ], [ chi[ PowerMap( mgmt, 3, i ) ],
    >                                  chi[ PowerMap( mgmt, 5, i ) ] ] );
    >     if v > 7 then
    >       v:= v - 15;
    >     fi;
    >     chi[i]:= AbsInt( v );
    >   od;
    gap> chi{ failpos };
    [ 2, 0, 0, 0, 1, 1, 2, 0, 0, 0, 1, 1 ]
\end{verbatim}


In order to assign the right labels to the preimages of those classes
of $C / \langle z \rangle$ that split into two classes of $C$,
we use congruences w.~r.~t.~the third power map.
Note that we have chosen that \emph{positive}
values of the projective character
belong to the \emph{first} entry in each pair of class labels.

\begin{verbatim}
    gap> split:= Filtered( failpos, x -> chi[x] <> 0 );
    [ 343, 347, 348, 383, 387, 388 ]
    gap> nonsplit:= Difference( failpos, split );
    [ 344, 345, 346, 384, 385, 386 ]
    gap> pow:= PowerMap( mgmt, 3 ){ split };
    [ 138, 136, 136, 263, 261, 261 ]
    gap> chi{ split };
    [ 2, 1, 1, 2, 1, 1 ]
    gap> chi{ pow };
    [ 8, 2, 2, 4, 2, 2 ]
    gap> cent_table.fusionlabels{ split };
    [ [ "15A", "30D" ], [ "15B", "30GH" ], [ "15B", "30GH" ], [ "30A", "30E" ], 
      [ "30GH", "30F" ], [ "30GH", "30F" ] ]
    gap> cent_table.fusionlabels{ pow };
    [ [ "5A", "10B" ], [ "10D", "5B" ], [ "10D", "5B" ], [ "10A", "10E" ], 
      [ "10F", "10D" ], [ "10F", "10D" ] ]
\end{verbatim}

We see from the element orders that the first three pairs are already
sorted correctly.
The fourth pair must be swapped because \texttt{30A} cubes to \texttt{10A}
and the positive values of the projective character
are not congruent modulo $3$.
Similarly, the last two pairs need not be swapped because \texttt{30F}
cubes to \texttt{10F}.

\begin{verbatim}
    gap> First( powerinfo, l -> l[1] = "30A" );
    [ "30A", [ [ 2, "15A" ], [ 3, "10A" ], [ 5, "6B" ] ] ]
    gap> First( powerinfo, l -> l[1] = "30F" );
    [ "30F", [ [ 2, "15B" ], [ 3, "10F" ], [ 5, "6I" ] ] ]
    gap> cent_table.fusionlabels[ split[4] ]:= Permuted(
    >        cent_table.fusionlabels[ split[4] ], (1,2) );;
    gap> for i in nonsplit do
    >      Unbind( cent_table.fusionlabels[i][2] );
    >    od;
\end{verbatim}

Assuming that not more classes of $C / \langle z \rangle$ split
under the epimorphism from $C$,
we create the (preliminary) character table head of $C$.

\begin{verbatim}
    gap> for i in split do
    >      cent_table.preimages[i]:= [ i, i ];
    >    od;
    gap> factorfusion:= Concatenation( cent_table.preimages );;
    gap> cb2b:= rec( UnderlyingCharacteristic:= 0,
    >                Size:= 2 * Size( mgmt ),
    >                Identifier:= "C_B(2B)",
    >                SizesCentralizers:= [] );;
    gap> proj:= InverseMap( factorfusion );;
    gap> for i in [ 1 .. Length( proj ) ] do
    >      if IsInt( proj[i] ) then
    >        cb2b.SizesCentralizers[ proj[i] ]:= SizesCentralizers( mgmt )[i];
    >      else
    >        cb2b.SizesCentralizers{ proj[i] }:=
    >            SizesCentralizers( mgmt )[i] * [ 2, 2 ];
    >      fi;
    >    od;
\end{verbatim}

The element orders are given by the class labels.

\begin{verbatim}
    gap> fusionlabels:= Concatenation( cent_table.fusionlabels );;
    gap> cb2b.OrdersClassRepresentatives:= List(
    >        Concatenation( cent_table.fusionlabels ),
    >        x -> Int( Filtered( x, IsDigitChar ) ) );;
    gap> ConvertToCharacterTable( cb2b );;
\end{verbatim}

Next we turn the projective character \texttt{chi} of $C/\langle z \rangle$
into a character for $C$,
and inflate the characters of the factor group to $C$.

\begin{verbatim}
    gap> chi_c:= [];;
    gap> for i in [ 1 .. Length( chi ) ] do
    >      if chi[i] = 0 then
    >        chi_c[ proj[i] ]:= 0;
    >      else
    >        chi_c{ proj[i] }:= [ 1, -1 ] * chi[i];
    >      fi;
    >    od;
    gap> factirr:= List( Irr( mgmt ), x -> x{ factorfusion } );;
\end{verbatim}

Tensoring the faithful irreducible character with the $60$ irreducible
characters of the factor group $Co_2$ yields $60$ irreducible characters
of $C$, which are in fact all the missing irreducibles.
In particular, no further splitting of classes occurs.

\begin{verbatim}
    gap> factirr_co2:= Filtered( factirr,
    >                      x -> ClassPositionsOfKernel( x ) <> [ 1, 2 ] );;
    gap> Length( factirr_co2 );
    60
    gap> ten:= Tensored( factirr_co2, [ chi_c ] );;
    gap> Set( List( ten, x -> ScalarProduct( cb2b, x, x ) ) );
    [ 1 ]
    gap> irr:= Concatenation( factirr, ten );;
    gap> Size( cb2b ) = Sum( List( irr, x -> x[1]^2 ) );
    true
    gap> SetIrr( cb2b, List( irr, x -> Character( cb2b, x ) ) );
\end{verbatim}

The only missing information on the character table of $C$
is that on power maps.
We use the power maps of the table for $C / \langle z \rangle$
as approximations, and let the standard algorithms compute the
candidates for the maps of $C$;
for all primes, there is only one candidate.


\begin{verbatim}
    gap> powermaps:= ComputedPowerMaps( cb2b );
    [  ]
    gap> for p in Set( Factors( Size( cb2b ) ) ) do
    >      init:= CompositionMaps( InverseMap( factorfusion ),
    >          CompositionMaps( PowerMap( mgmt, p ), factorfusion ) );
    >      poss:= PossiblePowerMaps( cb2b, p, rec( powermap:= init ) );
    >      if Length( poss ) <> 1 then
    >        Error( Ordinal( p ), " power map is not unique" );
    >      fi;
    >      powermaps[p]:= poss[1];
    >    od;
\end{verbatim}

Finally, we compare the newly computed character table
with that from {\GAP}'s library.

\begin{verbatim}
    gap> libtbl:= CharacterTable( "BM2" );;
    gap> Set( Irr( libtbl ) ) = Set( Irr( cb2b ) );
    true
    gap> IsRecord( TransformingPermutationsCharacterTables( libtbl, cb2b ) );
    true
\end{verbatim}

\section{Conjugacy classes of $\B$ and their centralizer orders}%
\label{sect:classes}

In this section, we determine the conjugacy classes of $\B$
and their centralizer orders, using the fact that for each element $g$
(of prime order) in a group $G$, say,
the conjugacy classes of elements in $G$ that contain roots of $g$
are in bijection with the conjugacy classes in $N_G(\langle g \rangle)$
that contain roots of $g$,
and that this bijection respects centralizer orders.

%

We have the following information about the elements of odd prime
order in $\B$, and their normalizers, see~\cite[Section 6]{BMverify}.

\begin{itemize}
\item
    There are exactly two classes of element order $3$ in $\B$,
    \texttt{3A} with normalizer $S_3 \times Fi_{22}.2$ and
    \texttt{3B} with normalizer $3^{1+8}.2^{1+6}.U_4(2).2$.
    Both subgroups have been constructed explicitly in the certified copy
    of $\B$, and the character table of the second subgroup has been
    recomputed from a permutation representation --it is equivalent to the
    character table in {\GAP}'s library of character tables.
\item
    There are exactly two classes of element order $5$ in $\B$,
    \texttt{5A} with normalizer $5:4 \times HS$ and
    \texttt{5B} with normalizer $5^{1+4}.2^{1+4}.A_5.4$.
    According to~\cite[Section~6]{BMverify},
    the two subgroups have been constructed explicitly in the certified copy
    of $\B$, and the character table of the second subgroup has been
    recomputed from a permutation representation --it is equivalent to the
    character table in {\GAP}'s library of character tables.
\item
    There is exactly one class of element order $7$ in $\B$,
    with normalizer of order $2^9 \cdot 3^3 \cdot 5 \cdot 7^2$
    and contained in maximal subgroups of type $2.{}^2E_6(2).2$.
\item
    There is exactly one class of element order $11$ in $\B$,
    with normalizer of type $S_5 \times 11:10$.
\item
    There is exactly one class of element order $13$ in $\B$,
    with normalizer of type $S_4 \times 13:12$.
\item
    There is exactly one class of element order $17$ in $\B$,
    with normalizer of order $2^6 \cdot 17$.
\item
    There is exactly one class of element order $19$ in $\B$,
    with normalizer of order $2^2 \cdot 3^2 \cdot 19$.
\item
    There are exactly two classes of element order $23$ in $\B$,
    which are Galois conjugate and have normalizer $2 \times 23:11$.
\item
    There are exactly two classes of element order $31$ in $\B$,
    which are Galois conjugate and have normalizer $31:15$.
\item
    There are exactly two classes of element order $47$ in $\B$,
    which are Galois conjugate and have normalizer $47:23$.
\end{itemize}

We are going to create the character table head for $\B$.
For that, we collect the information about already identified classes
in a record; its component names are the class names,
the corresponding values are the element order and the centralizer order.
The following auxiliary function sets a value in this record
(and signals an error if the value constradicts the one stored in the
library character tableof $\B$).

\begin{verbatim}
    gap> libB:= CharacterTable( "B" );;
    gap> libBnames:= ClassNames( libB, "ATLAS" );;
    gap> Bclassinfo:= rec();;
    gap> SetCentralizerOrder:= function( classname, value )
    >        local pos;
    > 
    >        pos:= Position( libBnames, classname );
    >        if pos = fail then
    >          Print( "no class name '", classname, "'?" );
    >          return false;
    >        elif SizesCentralizers( libB )[ pos ] <> value then
    >          Error( "wrong centralizer order!" );
    >        fi;
    >        Bclassinfo.( classname ):=
    >            [ Int( Filtered( classname, IsDigitChar ) ), value ];
    >        return true;
    >    end;;
\end{verbatim}

The values mentioned above, for elements of prime order at least $11$,
are entered by hand.

\begin{verbatim}
    gap> SetCentralizerOrder( "1A", Size( libB ) );;
    gap> SetCentralizerOrder( "11A", 11 * 120 );;
    gap> SetCentralizerOrder( "13A", 13 * 24 );;
    gap> SetCentralizerOrder( "17A", 17 * 4 );;
    gap> SetCentralizerOrder( "19A", 19 * 2 );;
    gap> SetCentralizerOrder( "23A", 23 * 2 );;
    gap> SetCentralizerOrder( "23B", 23 * 2 );;
    gap> SetCentralizerOrder( "31A", 31 );;
    gap> SetCentralizerOrder( "31B", 31 );;
    gap> SetCentralizerOrder( "47A", 47 );;
    gap> SetCentralizerOrder( "47B", 47 );;
\end{verbatim}

We will need the information about how many root classes a given class
of $\B$ has,
based on the power map information for the labels
(the list \texttt{powerinfo}).

\begin{verbatim}
    gap> RootInfoFromLabels:= function( label )
    >     local found, exp, res, entry, pos, root, ord;
    > 
    >     found:= [ label ];
    >     exp:= [ Int( Filtered( label, IsDigitChar ) ) ];
    >     res:= rec( total:= [], labels:= [] );
    >     res.total[ exp[1] ]:= 1;
    >     res.labels[ exp[1] ]:= [ label ];
    >     for entry in powerinfo do
    >       for root in entry[2] do
    >         if not entry[1] in found then
    >           pos:= Position( found, root[2] );
    >           if pos <> fail then
    >             ord:= exp[ pos ] * root[1];
    >             Add( exp, ord );
    >             Add( found, entry[1] );
    >             if not IsBound( res.total[ ord ] ) then
    >               res.total[ ord ]:= 0;
    >               res.labels[ ord ]:= [];
    >             fi;
    >             res.total[ ord ]:= res.total[ ord ] + 1;
    >             Add( res.labels[ ord ], entry[1] );
    >           fi;
    >         fi;
    >       od;
    >     od;
    >     return res;
    > end;;
\end{verbatim}

Analogously, we will need this roots information for the normalizers
of certain elements, computed from the character tables.

\begin{verbatim}
    gap> RootInfoFromTable:= function( tbl, pos )
    >     local orders, posord, res, i, ord;
    > 
    >     orders:= OrdersClassRepresentatives( tbl );
    >     posord:= orders[ pos ];
    >     res:= rec( total:= [], classpos:= [] );
    >     for i in [ 1 .. NrConjugacyClasses( tbl ) ] do
    >       ord:= orders[i];
    >       if ord mod posord = 0 and
    >          PowerMap( tbl, ord / posord, i ) = pos then
    >         if not IsBound( res.total[ ord ] ) then
    >           res.total[ ord ]:= 0;
    >           res.classpos[ ord ]:= [];
    >         fi;
    >         res.total[ ord ]:= res.total[ ord ] + 1;
    >         Add( res.classpos[ ord ], i );
    >       fi;
    >     od;
    >     return res;
    > end;;
\end{verbatim}

We compute the roots info for the four involution normalizers,
from the available character tables.

\begin{verbatim}
    gap> norm2A:= CharacterTable( "2.2E6(2).2" );;
    gap> pos2A:= ClassPositionsOfCentre( norm2A );
    [ 1, 2 ]
    gap> rootinfo2A_t:= RootInfoFromTable( norm2A, pos2A[2] );;
    gap> pos2B:= ClassPositionsOfCentre( cb2b );
    [ 1, 2 ]
    gap> rootinfo2B_t:= RootInfoFromTable( cb2b, pos2B[2] );;
    gap> norm2C:= CharacterTableOfIndexTwoSubdirectProduct(
    >      CharacterTable( "F4(2)" ), CharacterTable( "F4(2).2" ),
    >      CharacterTable( "2^2" ), CharacterTable( "D8" ), "norm2C" );;
    gap> pos2C:= ClassPositionsOfCentre( norm2C.table );
    [ 1, 2 ]
    gap> rootinfo2C_t:= RootInfoFromTable( norm2C.table, pos2C[2] );;
    gap> sup_norm2D:= CharacterTable( "BM4" );
    CharacterTable( "2^(9+16).S8(2)" )
    gap> pos2D:= Positions( SizesCentralizers( sup_norm2D ), 11689182992793600 );
    [ 4 ]
    gap> rootinfo2D_t:= RootInfoFromTable( sup_norm2D, pos2D[1] );;
    gap> rootinfo2A_l:= RootInfoFromLabels( "2A" );;
    gap> rootinfo2B_l:= RootInfoFromLabels( "2B" );;
    gap> rootinfo2C_l:= RootInfoFromLabels( "2C" );;
    gap> rootinfo2D_l:= RootInfoFromLabels( "2D" );;
\end{verbatim}

We have to justify that the labels \texttt{"2A"}, \texttt{"2B"},
\texttt{"2C"}, \texttt{"2D"} for class representatives fit to the
class names \texttt{2A}, \texttt{2B}, \texttt{2C}, \texttt{2D}
which are used in \cite{Str76b}.
For that, we check for the availability of roots of certain orders.

\begin{verbatim}
    gap> List( [ rootinfo2A_t, rootinfo2B_t, rootinfo2C_t, rootinfo2D_t ],
    >          r -> Filtered( [ 52, 56, 70 ], i -> IsBound( r.total[i] ) ) );
    [ [ 70 ], [ 56 ], [ 52 ], [  ] ]
    gap> List( [ rootinfo2A_l, rootinfo2B_l, rootinfo2C_l, rootinfo2D_l ],
    >          r -> Filtered( [ 52, 56, 70 ], i -> IsBound( r.total[i] ) ) );
    [ [ 70 ], [ 56 ], [ 52 ], [  ] ]
\end{verbatim}

Now we can start to identify the root classes of the normalizers
with the class labels.

\begin{verbatim}
    gap> List( [ rootinfo2A_t, rootinfo2B_t, rootinfo2C_t, rootinfo2D_t ],
    >          r -> Sum( Compacted( r.total ) ) );
    [ 20, 83, 14, 40 ]
    gap> List( [ rootinfo2A_l, rootinfo2B_l, rootinfo2C_l, rootinfo2D_l ],
    >          r -> Sum( Compacted( r.total ) ) );
    [ 19, 77, 14, 40 ]
\end{verbatim}

We see that some of the class labels for $\B$ belong to more than one
conjugacy class.  Exactly one such case occurs for a root class of
\texttt{2A}, the label \texttt{"34BC"} belongs to two Galois conjugate
classes \texttt{34B}, \texttt{34C}.
(This was already clear from the discussion of classes of elements
of order $17$.)

\begin{verbatim}
    gap> rootinfo2A_t.total[34];
    2
    gap> rootinfo2A_l.total[34];
    1
    gap> rootinfo2A_l.labels[34];
    [ "34BC" ]
    gap> pos34:= rootinfo2A_t.classpos[34];
    [ 133, 135 ]
    gap> PowerMap( norm2A, 3 ){ pos34 };
    [ 135, 133 ]
\end{verbatim}

The other six cases where we have to adjust the labels occur for roots of
\texttt{2B}.

\begin{verbatim}
    gap> rootinfo2B_t.total{ [ 16, 30, 32, 46, 56 ] };
    [ 6, 3, 4, 2, 2 ]
    gap> rootinfo2B_l.total{ [ 16, 30, 32, 46, 56 ] };
    [ 5, 2, 2, 1, 1 ]
\end{verbatim}

For the classes of element order different from $16$,
we have to replace a class label by a pair of labels
which belong to Galois conjugate classes:

\begin{verbatim}
    gap> rootinfo2B_l.labels[30];
    [ "30GH", "30D" ]
    gap> pos30:= rootinfo2B_t.classpos[30];
    [ 388, 393, 395 ]
    gap> PowerMap( cb2b, 7 ){ pos30 };
    [ 388, 395, 393 ]
    gap> List( pos30,
    >          i -> Number( PowerMap( cb2b, 2 ), x -> x = i ) );
    [ 2, 0, 0 ]
    gap> List( rootinfo2B_l.labels[30],
    >          l -> Length( RootInfoFromLabels( l ).total ) );
    [ 30, 60 ]
\end{verbatim}

We have to replace the label \texttt{"30GH"} (for which there are no square
roots) into two labels, which we call \texttt{"30G"} and \texttt{"30H"}.

\begin{verbatim}
    gap> rootinfo2B_l.labels[32];
    [ "32AB", "32CD" ]
    gap> pos32:= rootinfo2B_t.classpos[32];
    [ 399, 400, 403, 404 ]
    gap> PowerMap( cb2b, 5 ){ pos32 };
    [ 400, 399, 404, 403 ]
\end{verbatim}

We replace \texttt{"32AB"} by \texttt{"32A"} and \texttt{"32B"},
and replace \texttt{"32GH"} by \texttt{"32G"} and \texttt{"32H"}.

\begin{verbatim}
    gap> rootinfo2B_l.labels[46];
    [ "46AB" ]
\end{verbatim}

It was already clear from the presence of two Galois conjugate classes
of element order $23$ that there must be two Galois conjugate classes
of element order $46$;
we replace \texttt{"46AB"} by \texttt{"46A"} and \texttt{"46B"}.

\begin{verbatim}
    gap> rootinfo2B_l.labels[56];
    [ "56AB" ]
    gap> pos56:= rootinfo2B_t.classpos[56];
    [ 437, 438 ]
    gap> PowerMap( cb2b, 3 ){ pos56 };
    [ 437, 438 ]
    gap> PowerMap( cb2b, 5 ){ pos56 };
    [ 438, 437 ]
\end{verbatim}

We replace \texttt{"56AB"} by \texttt{"56A"} and \texttt{"56B"}.

Now one case of element order $16$ is left.
There are two classes of element order $8$ for which we
need three classes of square roots.
We have three square roots for the label \texttt{"8D"}
but only two for \texttt{"8H"},
which are \texttt{"16C"} and \texttt{"16DF"}.
Since both of them have square roots,
we conclude that we have to introduce one new label for a class
of element order $16$ that has no square roots.
We call this new label \texttt{"16F"},
and replace \texttt{"16DF"} by \texttt{"16D"}.

\begin{verbatim}
    gap> pos:= rootinfo2B_t.classpos[16];
    [ 233, 234, 235, 251, 252, 258 ]
    gap> PowerMap( cb2b, 2 ){ pos };
    [ 69, 69, 69, 104, 104, 104 ]
    gap> List( rootinfo2B_t.classpos[16],
    >          i -> Number( PowerMap( cb2b, 2 ), x -> x = i ) );
    [ 0, 0, 0, 2, 0, 2 ]
    gap> rootinfo2B_l.labels[16];
    [ "16B", "16C", "16A", "16E", "16DF" ]
    gap> Filtered( powerinfo, l -> l[1] in rootinfo2B_l.labels[16] );
    [ [ "16B", [ [ 2, "8D" ] ] ], [ "16C", [ [ 2, "8H" ] ] ], 
      [ "16A", [ [ 2, "8D" ] ] ], [ "16E", [ [ 2, "8D" ] ] ], 
      [ "16DF", [ [ 2, "8H" ] ] ] ]
    gap> List( rootinfo2B_l.labels[16],
    >          l -> IsBound( RootInfoFromLabels( l ).total[32] ) );
    [ false, true, false, false, true ]
\end{verbatim}


We record the splittings in the list \texttt{powerinfo}.
First we replace the entry for each splitting label by two entries.

\begin{verbatim}
    gap> tosplit:= List( Filtered( powerinfo,
    >                              x -> Number( x[1], IsAlphaChar ) > 1 ),
    >                    x -> x[1] );
    [ "16DF", "23AB", "30GH", "31AB", "32AB", "32CD", "34BC", "46AB", "47AB", 
      "56AB" ]
    gap> for l in tosplit do
    >      ord:= Filtered( l, IsDigitChar );
    >      new:= List( Filtered( l, IsAlphaChar ),
    >                  c -> Concatenation( ord, [ c ] ) );
    >      pos:= PositionProperty( powerinfo, x -> x[1] = l );
    >      Append( powerinfo, List( new, x -> [ x,
    >                  StructuralCopy( powerinfo[ pos ][2] ) ] ) );
    >      Unbind( powerinfo[ pos ] );
    >    od;
    gap> powerinfo:= Compacted( powerinfo );;
    gap> SortParallel(
    >      List( powerinfo,
    >            x -> [ Int( Filtered( x[1], IsDigitChar ) ), x[1] ] ),
    >      powerinfo );
\end{verbatim}

The splitting classes occur as powers only in the follwoing cases:
We may choose \texttt{"23A"} as the square of \texttt{"46A"},
and \texttt{"23B"} as the square of \texttt{"46B"}.
And we have already said that the square of both \texttt{"32C"}
and \texttt{"32D"} shall be \texttt{"16D"}.
We adjust these cases by hand.

\begin{verbatim}
    gap> Filtered( powerinfo, x -> ForAny( x[2], p -> p[2] in tosplit ) );
    [ [ "32C", [ [ 2, "16DF" ] ] ], [ "32D", [ [ 2, "16DF" ] ] ], 
      [ "46A", [ [ 2, "23AB" ], [ 23, "2B" ] ] ], 
      [ "46B", [ [ 2, "23AB" ], [ 23, "2B" ] ] ] ]
    gap> entry:= First( powerinfo, x -> x[1] = "32C" );;
    gap> entry[2][1][2]:= "16D";;
    gap> entry:= First( powerinfo, x -> x[1] = "32D" );;
    gap> entry[2][1][2]:= "16D";;
    gap> entry:= First( powerinfo, x -> x[1] = "46A" );;
    gap> entry[2][1][2]:= "23A";;
    gap> entry:= First( powerinfo, x -> x[1] = "46B" );;
    gap> entry[2][1][2]:= "23B";;
\end{verbatim}

Now the numbers of roots for the four involution labels coincide
with the corresponding numbers of root classes in the normalizers.
We know that $\B$ has $184$ conjugacy classes.

\begin{verbatim}
    gap> Length( powerinfo );
    184
    gap> Bnames:= List( powerinfo, x -> x[1] );;
\end{verbatim}

Next we verify that the labels for elements of odd order
describe already the conjugacy classes of elements of odd order.
For that, it is sufficient to check the normalizers of
\texttt{"3A"}, \texttt{"3B"}, \texttt{"5A"}, \texttt{"5B"}.
Note that \texttt{"3A"} has roots of order $66$ and
\texttt{"3B"} has no such roots,
and \texttt{"5A"} has roots of order $70$ and \texttt{"5B"} has not;
this means that the names of the labels coincide with the class names.

\begin{verbatim}
    gap> n3a:= CharacterTable( "S3xFi22.2" );;
    gap> pos:= ClassPositionsOfPCore( n3a, 3 );
    [ 1, 113 ]
    gap> rootinfo3A_t:= RootInfoFromTable( n3a, pos[2] );;
    gap> rootinfo3A_l:= RootInfoFromLabels( "3A" );;
    gap> n3b:= CharacterTable( "3^(1+8).2^(1+6).U4(2).2" );;
    gap> pos:= ClassPositionsOfPCore( n3b, 3 );
    [ 1 .. 4 ]
    gap> Filtered( ClassPositionsOfNormalSubgroups( n3b ),
    >              n -> IsSubset( pos, n ) );
    [ [ 1 ], [ 1, 2 ], [ 1, 2, 3, 4 ] ]
    gap> rootinfo3B_t:= RootInfoFromTable( n3b, pos[2] );;
    gap> rootinfo3B_l:= RootInfoFromLabels( "3B" );;
    gap> IsBound( rootinfo3A_t.total[66] );
    true
    gap> IsBound( rootinfo3B_t.total[66] );
    false
    gap> rootinfo3A_t.total = rootinfo3A_l.total;
    true
    gap> rootinfo3B_t.total = rootinfo3B_l.total;
    true
    gap> n5a:= CharacterTable( "5:4xHS.2" );
    CharacterTable( "5:4xHS.2" )
    gap> pos:= ClassPositionsOfPCore( n5a, 5 );
    [ 1, 40 ]
    gap> rootinfo5A_t:= RootInfoFromTable( n5a, pos[2] );;
    gap> rootinfo5A_l:= RootInfoFromLabels( "5A" );;
    gap> n5b:= CharacterTable( "5^(1+4).2^(1+4).A5.4" );
    CharacterTable( "5^(1+4).2^(1+4).A5.4" )
    gap> pos:= ClassPositionsOfPCore( n5b, 5 );
    [ 1 .. 4 ]
    gap> Filtered( ClassPositionsOfNormalSubgroups( n5b ),
    >              n -> IsSubset( pos, n ) );
    [ [ 1 ], [ 1, 2 ], [ 1, 2, 3, 4 ] ]
    gap> rootinfo5B_t:= RootInfoFromTable( n5b, pos[2] );;
    gap> rootinfo5B_l:= RootInfoFromLabels( "5B" );;
    gap> IsBound( rootinfo5A_t.total[70] );
    true
    gap> IsBound( rootinfo5B_t.total[70] );
    false
    gap> rootinfo5A_t.total = rootinfo5A_l.total;
    true
    gap> rootinfo5B_t.total = rootinfo5B_l.total;
    true
\end{verbatim}

It will be useful to provide the \texttt{powerinfo} information
in a record.

\begin{verbatim}
    gap> powerinforec:= rec();;
    gap> for entry in powerinfo do
    >      powerinforec.( entry[1] ):= entry[2];
    >    od;
\end{verbatim}

We recompute the roots information, according to the changed data.

\begin{verbatim}
    gap> rootinfo2A_l:= RootInfoFromLabels( "2A" );;
    gap> rootinfo2B_l:= RootInfoFromLabels( "2B" );;
    gap> rootinfo2C_l:= RootInfoFromLabels( "2C" );;
    gap> rootinfo2D_l:= RootInfoFromLabels( "2D" );;
    gap> rootinfo2A_t.total = rootinfo2A_l.total;
    true
    gap> rootinfo2B_t.total = rootinfo2B_l.total;
    true
    gap> rootinfo2C_t.total = rootinfo2C_l.total;
    true
    gap> rootinfo2D_t.total = rootinfo2D_l.total;
    true
\end{verbatim}

We try to identify the classes with labels.
The numbers of classes and labels fit together,
now we compute the bijection.
The following function identifies classes which are determined either
already by the element order or as a power of an identified class
or as a unique root of an identified class.

\begin{verbatim}
    gap> IdentifyCentralizerOrders:= function( normtbl, rl, rt )
    >     local n, identified, found, i, unknown, class, d, linfo, p, e, cand,
    >           imgs, im, pos, powerlabel, dd, cent;
    >     n:= First( [ 1 .. Length( rl ) ], i -> IsBound( rl[i] ) );
    >     identified:= [ [], [] ];
    >     found:= true;
    >     while found do
    >       found:= false;
    >       for i in [ 1 .. Length( rl ) ] do
    >       if IsBound( rl[i] ) then
    >         unknown:= Difference( rl[i], identified[1] );
    >         if Length( unknown ) = 1 then
    >           # Identify the class.
    >           class:= Difference( rt[i], identified[2] )[1];
    >           Add( identified[1], unknown[1] );
    >           Add( identified[2], class );
    >           found:= true;
    >           # Identify the admissible powers.
    >           for d in Difference( DivisorsInt( i / n ), [ 1 ] ) do
    >             linfo:= powerinforec.( unknown[1] );
    >             for p in Factors( d ) do
    >               e:= First( linfo, x -> x[1] = p );
    >               linfo:= powerinforec.( e[2] );
    >             od;
    >             if not e[2] in identified[1] then
    >               Add( identified[1], e[2] );
    >               Add( identified[2], PowerMap( normtbl, d, class ) );
    >               found:= true;
    >             fi;
    >           od;
    >         else
    >           # Try to identify roots whose powers are identified.
    >           for d in Difference( DivisorsInt( i / n ), [ 1 ] ) do
    >             cand:= Difference( rt[i], identified[2] );
    >             imgs:= PowerMap( normtbl, d ){ cand };
    >             for im in Intersection( imgs, identified[2] ) do
    >               pos:= Positions( imgs, im );
    >               if Length( pos ) = 1 then
    >                 class:= cand[ pos[1] ];
    >                 powerlabel:= identified[1][
    >                                  Position( identified[2], im ) ];
    >                 # Find the labels of the 'd'-th powers of 'unknown'.
    >                 linfo:= List( unknown, l -> powerinforec.( l ) );
    >                 for p in Factors( d ) do
    >                   e:= List( linfo, ll -> First( ll, x -> x[1] = p ) );
    >                   linfo:= List( e, ee -> powerinforec.( ee[2] ) );
    >                 od;
    >                 linfo:= List( e, x -> x[2] );
    >                 pos:= Position( linfo, powerlabel );
    >                 Add( identified[1], unknown[ pos ] );
    >                 Add( identified[2], class );
    >                 # Identify the admissible powers.
    >                 for dd in Difference( DivisorsInt( i / n ), [ 1 ] ) do
    >                   linfo:= powerinforec.( unknown[ pos ] );
    >                   for p in Factors( dd ) do
    >                     e:= First( linfo, x -> x[1] = p );
    >                     linfo:= powerinforec.( e[2] );
    >                   od;
    >                   if not e[2] in identified[1] then
    >                     Add( identified[1], e[2] );
    >                     Add( identified[2], PowerMap( normtbl, dd, class ) );
    >                     found:= true;
    >                   fi;
    >                 od;
    >                 found:= true;
    >                 break; # since we have to update 'unknown'
    >               fi;
    >             od;
    >             if found then
    >               break; # since we have to update 'unknown'
    >             fi;
    >           od;
    >         fi;
    >       fi;
    >     od;
    >   od;
    >   # Where the centralizer order is unique, set it.
    >   for i in [ 1 .. Length( rl ) ] do
    >     if IsBound( rl[i] ) then
    >       cand:= Difference( rt[i], identified[2] );
    >       cent:= Set( SizesCentralizers( normtbl ){ cand } );
    >       if Length( cent ) = 1 then
    >         Append( identified[1], Difference( rl[i], identified[1] ) );
    >         Append( identified[2], cand );
    >       fi;
    >     fi;
    >   od;
    >   # Set the centralizer orders.
    >   for i in [ 1 .. Length( identified[1] ) ] do
    >     if not IsBound( Bclassinfo.( identified[1][i] ) ) then
    >       Print( "#I  identify ", identified[1][i], "\n" );
    >     fi;
    >     SetCentralizerOrder( identified[1][i],
    >         SizesCentralizers( normtbl )[ identified[2][i] ] );
    >   od;
    >   # Return the information about unidentified classes.
    >   return [ Difference( Concatenation( Compacted( rl ) ), identified[1] ),
    >            Difference( Concatenation( Compacted( rt ) ), identified[2] ) ];
    > end;;
\end{verbatim}

We try the function with the four involution normalizers.
In the case of \texttt{cb2b},
we are better off since we know most of the class fusion to $\B$.
Thus we use also the information about the labels that belong to roots of
\texttt{"2B"} elements where this is available.

\begin{verbatim}
    gap> IdentifyCentralizerOrders( norm2A,
    >        rootinfo2A_l.labels, rootinfo2A_t.classpos );
    #I  identify 2A
    #I  identify 10A
    #I  identify 22A
    #I  identify 26B
    #I  identify 38A
    #I  identify 66A
    #I  identify 6A
    #I  identify 70A
    #I  identify 14A
    #I  identify 14B
    #I  identify 42A
    #I  identify 6B
    #I  identify 6D
    #I  identify 42B
    #I  identify 30B
    #I  identify 30A
    #I  identify 34B
    #I  identify 34C
    [ [ "18A", "18B" ], [ 74, 76 ] ]
    gap> for i in Union( rootinfo2B_t.classpos ) do
    >      if IsBound( powerinforec.( fusionlabels[i] ) ) then
    >        SetCentralizerOrder( fusionlabels[i],
    >            SizesCentralizers( cb2b )[i] );
    >      fi;
    >    od;
    gap> IdentifyCentralizerOrders( cb2b,
    >        rootinfo2B_l.labels, rootinfo2B_t.classpos );
    #I  identify 30G
    #I  identify 30H
    #I  identify 32A
    #I  identify 32B
    #I  identify 32C
    #I  identify 32D
    #I  identify 46A
    #I  identify 46B
    #I  identify 56A
    #I  identify 56B
    [ [ "12A", "12D", "12G", "12I", "12L", "12M", "12O", "16A", "16B", "16C", 
          "16D", "16E", "16F", "20B", "20C", "20F", "20I", "24A", "24B", "24C", 
          "24D", "24E", "24F", "24G", "24K", "24M", "40A", "40B", "40C", "40D", 
          "4F", "4G", "8A", "8B", "8C", "8D", "8E", "8F", "8H", "8I", "8L" ], 
      [ 28, 42, 47, 48, 49, 61, 63, 68, 69, 103, 104, 116, 145, 151, 158, 163, 
          169, 187, 199, 215, 217, 218, 219, 220, 233, 234, 235, 251, 252, 258, 
          292, 308, 309, 310, 311, 320, 332, 333, 344, 357, 420 ] ]
    gap> IdentifyCentralizerOrders( norm2C.table,
    >        rootinfo2C_l.labels, rootinfo2C_t.classpos );
    #I  identify 2C
    #I  identify 4I
    #I  identify 10C
    #I  identify 12T
    #I  identify 6K
    #I  identify 14C
    #I  identify 18F
    #I  identify 20H
    #I  identify 26A
    #I  identify 30C
    #I  identify 6F
    #I  identify 34A
    #I  identify 42C
    #I  identify 52A
    [ [  ], [  ] ]
    gap> IdentifyCentralizerOrders( sup_norm2D,
    >        rootinfo2D_l.labels, rootinfo2D_t.classpos );
    #I  identify 2D
    #I  identify 14E
    #I  identify 28E
    #I  identify 4E
    #I  identify 36C
    #I  identify 18E
    #I  identify 12N
    #I  identify 6J
    #I  identify 40E
    #I  identify 20G
    #I  identify 10F
    #I  identify 8G
    #I  identify 60C
    #I  identify 30F
    #I  identify 12F
    #I  identify 6I
    #I  identify 8J
    #I  identify 10E
    #I  identify 12S
    #I  identify 4H
    #I  identify 18D
    #I  identify 20J
    #I  identify 4J
    #I  identify 24H
    #I  identify 30E
    #I  identify 6E
    #I  identify 6H
    #I  identify 8N
    #I  identify 12J
    #I  identify 12P
    #I  identify 12Q
    #I  identify 24L
    #I  identify 8K
    #I  identify 8M
    #I  identify 12R
    #I  identify 16G
    #I  identify 24N
    #I  identify 16H
    #I  identify 24I
    #I  identify 24J
    [ [  ], [  ] ]
\end{verbatim}

\begin{verbatim}
    gap> IdentifyCentralizerOrders( n3a,
    >        rootinfo3A_l.labels, rootinfo3A_t.classpos );;
    #I  identify 3A
    #I  identify 15A
    #I  identify 21A
    #I  identify 33A
    #I  identify 39A
    gap> IdentifyCentralizerOrders( n3b,
    >        rootinfo3B_l.labels, rootinfo3B_t.classpos );;
    #I  identify 3B
    #I  identify 15B
    #I  identify 27A
    #I  identify 9A
    #I  identify 9B
    gap> IdentifyCentralizerOrders( n5a,
    >        rootinfo5A_l.labels, rootinfo5A_t.classpos );;
    #I  identify 5A
    #I  identify 35A
    #I  identify 55A
    gap> IdentifyCentralizerOrders( n5b,
    >        rootinfo5B_l.labels, rootinfo5B_t.classpos );;
    #I  identify 5B
    #I  identify 25A
\end{verbatim}

Let us see which classes of $\B$ are not identified yet.

\begin{verbatim}
    gap> Difference( RecNames( powerinforec ), RecNames( Bclassinfo ) );
    [ "16D", "16F", "18A", "18B", "7A" ]
\end{verbatim}

For simplicity, we set the centralizer order for \texttt{"7A"} by hand.

\begin{verbatim}
    gap> SetCentralizerOrder( "7A", 2^8 * 3^2 * 5 * 7^2 );;
\end{verbatim}

The classes \texttt{18A} and \texttt{"18B"} are roots of \texttt{"2A"}.
They have no roots and the same power maps.
The centralizer orders of the two classes in the \texttt{"2A"} centralizer
are $1296 = 2^4 \cdot 3^4$ and $648 = 2^3 \cdot 3^4$, respectively.
We know a class in \texttt{cb2b} that fuses into \texttt{18A} and
has centralizer order $144 = 2^4 \cdot 3^2$ in \texttt{cb2b}.
Thus \texttt{18A} must have centralizer order $1296$ in $\B$.

\begin{verbatim}
    gap> SetCentralizerOrder( "18A", 1296 );;
    gap> SetCentralizerOrder( "18B", 648 );;
\end{verbatim}

The cases \texttt{"16D"} and \texttt{"16F"} will be handled below.

Later we will need the class fusion from $C_B(\texttt{2B})$ to $\B$,
and we know it already as far as the class invariants reach.
We compute part of the missing information.

\begin{verbatim}
    gap> diff:= Difference( fusionlabels, RecNames( powerinforec ) );
    [ "16DF", "23AB", "30GH", "32AB", "32CD", "46AB", "56AB" ]
    gap> List( diff, x -> Positions( fusionlabels, x ) );
    [ [ 251, 252, 274, 281, 396 ], [ 421, 423 ], [ 393, 395, 445, 447 ], 
      [ 399, 400 ], [ 403, 404 ], [ 422, 424 ], [ 437, 438 ] ]
    gap> Length( fusionlabels );
    448
    gap> NrConjugacyClasses( cb2b );
    448
\end{verbatim}

We are free to choose the images of the class fusion for the elements of
order $23$ (which then determines the classes of element order $46$),
$32$, and $56$, since the question is about independent pairs
of Galois conjugate classes.

\begin{verbatim}
    gap> fusionlabels[421]:= "23A";;
    gap> fusionlabels[423]:= "23B";;
    gap> fusionlabels[399]:= "32A";;
    gap> fusionlabels[400]:= "32B";;
    gap> fusionlabels[403]:= "32C";;
    gap> fusionlabels[404]:= "32D";;
    gap> fusionlabels[437]:= "56A";;
    gap> fusionlabels[438]:= "56B";;
    gap> pos46:= Positions( OrdersClassRepresentatives( cb2b ), 46 );
    [ 422, 424 ]
    gap> PowerMap( cb2b, 2 ){ pos46 };
    [ 421, 423 ]
    gap> fusionlabels[422]:= "46A";;
    gap> fusionlabels[424]:= "46B";;
\end{verbatim}

Thus we are left with the question about the fusion to the classes with the
labels \texttt{"16D"}, \texttt{"16F"}, \texttt{"30G"}, and \texttt{"30H"}.

For the order $30$ elements, we may choose images for \emph{one} pair of
Galois conjugate classes, and later try to distinguish the two possibilities
for the other pair, for example via induced characters.

\begin{verbatim}
    gap> fusionlabels[393]:= "30G";;
    gap> fusionlabels[395]:= "30H";;
\end{verbatim}


Two classes of order $16$ elements with fusion label \texttt{"16DF"}
are roots of the central involution of \texttt{cb2b},
and we can distinguish them by the fact that
\texttt{"16D"} has square roots whereas \texttt{"16F"} has not.
For the other three classes, we are left with two possibilities.

\begin{verbatim}
    gap> fusionlabels[251]:= "16D";;
    gap> fusionlabels[252]:= "16F";;
    gap> SetCentralizerOrder( "16D", SizesCentralizers( cb2b )[251] );;
    gap> SetCentralizerOrder( "16F", SizesCentralizers( cb2b )[252] );;
\end{verbatim}

This means that we have currently $2^5$ candidates for the class fusion
from \texttt{cb2b} to $\B$.

\begin{verbatim}
    gap> cb2bfusb:= List( fusionlabels, l -> Position( Bnames, l ) );;
    gap> Positions( cb2bfusb, fail );
    [ 274, 281, 396, 445, 447 ]
    gap> pos:= Positions( fusionlabels, "16DF" );
    [ 274, 281, 396 ]
    gap> for i in pos do
    >      cb2bfusb[i]:= [ 86, 88 ];
    >    od;
    gap> pos:= Positions( fusionlabels, "30GH" );
    [ 445, 447 ]
    gap> for i in pos do
    >      cb2bfusb[i]:= [ 143, 144 ];
    >    od;
\end{verbatim}

Before we compute the irreducible characters of $\B$,
we create the character table head for $\B$.

\begin{verbatim}
    gap> bhead:= rec( UnderlyingCharacteristic:= 0,
    >                 Size:= Bclassinfo.( "1A" )[2],
    >                 Identifier:= "Bnew" );;
    gap> bhead.SizesCentralizers:= List( Bnames, x -> Bclassinfo.( x )[2] );;
    gap> bhead.OrdersClassRepresentatives:= List( Bnames,
    >        x -> Bclassinfo.( x )[1] );;
    gap> bhead.ComputedPowerMaps:= [];;
    gap> galoisinfo:= rec(
    >     classes:= [ "23A", "23B", "30G", "30H", "31A", "31B", "32A", "32B",
    >                 "32C", "32D", "34B", "34C", "46A", "46B", "47A", "47B",
    >                 "56A", "56B" ],
    >     partners:= [ "23B", "23A", "30H", "30G", "31B", "31A", "32B", "32A",
    >                  "32D", "32C", "34C", "34B", "46B", "46A", "47B", "47A",
    >                  "56B", "56A" ],
    >     rootsof:= [ -23, -23, -15, -15, -31, -31, 2, 2,
    >                 -2, -2, 17, 17, -23, -23, -47, -47,
    >                 7, 7 ] );;
    gap> galoisinfo.irrats:= List( galoisinfo.rootsof, Sqrt );;
    gap> for p in Filtered( [ 1 .. Maximum( bhead.OrdersClassRepresentatives ) ],
    >                       IsPrimeInt ) do
    >      map:= [ 1 ];
    >      for i in [ 2 .. Length( Bnames ) ] do
    >        if bhead.OrdersClassRepresentatives[i] = p then
    >          map[i]:= 1;
    >        elif bhead.OrdersClassRepresentatives[i] mod p = 0 then
    >          # The 'p'-th power has smaller order, we know the image class.
    >          info:= First( powerinforec.( Bnames[i] ), pair -> pair[1] = p );
    >          map[i]:= Position( Bnames, info[2] );
    >        else
    >          # The 'p'-th power is a Galois conjugate.
    >          pos:= Position( galoisinfo.classes, Bnames[i] );
    >          if pos = fail then
    >            # The 'i'-th class is rational.
    >            map[i]:= i;
    >          else
    >            # Determine whether the pair gets swapped.
    >            irrat:= galoisinfo.irrats[ pos ];
    >            if GaloisCyc( irrat, p ) <> irrat then
    >              map[i]:= Position( Bnames, galoisinfo.partners[ pos ] );
    >            else
    >              map[i]:= i;
    >            fi;
    >          fi;
    >        fi;
    >      od;
    >      bhead.ComputedPowerMaps[p]:= map;
    >    od;
    gap> ConvertToCharacterTable( bhead );;
\end{verbatim}

Check the library table of $\B$ against the character table head.

\begin{verbatim}
    gap> b:= CharacterTable( "B" );;
    gap> for p in Filtered( [ 2 ..
    >                 Maximum( OrdersClassRepresentatives( b ) ) ],
    >                 IsPrimeInt ) do
    >      PowerMap( b, p );
    >    od;
    gap> ComputedPowerMaps( bhead ) = ComputedPowerMaps( b );
    true
\end{verbatim}

Determine the missing pieces of the class fusion
from the \texttt{2B} centralizer to $\B$.

\begin{verbatim}
    gap> maps:= ContainedMaps( cb2bfusb );;
    gap> Length( maps );
    32
    gap> good:= [];;
    gap> for map in maps do
    >      ind:= InducedClassFunctionsByFusionMap( cb2b, b, Irr( cb2b ), map );
    >      if ForAll( ind, x -> IsInt( ScalarProduct( b, x, x ) ) ) then
    >        Add( good, map );
    >      fi;
    >    od;
    gap> Length( good );
    1
    gap> b2b:= CharacterTable( "BM2" );;
    gap> good[1] = GetFusionMap( b2b, b );
    true
\end{verbatim}

In our situation (where the classes of the subgroup have been
identified via the class invariants for $\B$,
and where we have made appropriate choices for the pairs of
Galois conjugate classes),
the class fusion is unique.

The newly computed fusion coincides with the one that is stored
on the {\GAP} library table.

\section{The irreducible characters of $\B$}\label{sect:irreducibles}

We assume the following information about the Baby Monster group $\B$.

\begin{itemize}
\item
    The conjugacy class lengths, the element orders, and the power maps
    (for all primes up to the maximal element order in $\B$) are known
    and coincide with the information that is shown in~\cite{CCN85}.
\item
    The group $\B$ contains subgroups of the structures
    $2.{}^2E_6(2).2$, $Fi_{23}$, and $HN.2$.
    The ordinary character tables of these groups have been verified
    (see~\cite{BMO17}) and thus may be used in our computations.
\item
    The character table of the \texttt{2B} centralizer in $\B$ is known.
    Also the class fusion from this table to the table head of $\B$
    is known by the construction of this character table
    in Section~\ref{table_c2b}.
\end{itemize}

For the sake of simplicity, we start with the {\ATLAS} table of $\B$
and store the power maps up to the maximal element order
(needed for inducing from cyclic subgroups).

\begin{verbatim}
    gap> b:= CharacterTable( "B" );;
    gap> for p in Filtered( [ 2 ..
    >                 Maximum( OrdersClassRepresentatives( b ) ) ],
    >                 IsPrimeInt ) do
    >      PowerMap( b, p );
    >    od;
\end{verbatim}

In order to make sure that the irreducible characters that are stored
on the table are not silently used inside some computations,
we delete them from the character table.

\begin{verbatim}
    gap> irr_atlas:= Irr( b );;
    gap> ResetFilterObj( b, HasIrr );
\end{verbatim}

Now we compute candidates for the class fusions from the subgroups
which we are allowed to use.
For that, we write a small {\GAP} program.
The input parameters are the character tables of the subgroup and $\B$,
a list of characters of $\B$, and perhaps a first approximation of the
class fusion in question.

\begin{verbatim}
    gap> tryFusion:= function( s, b, ind, initmap )
    >      local i, sfusb, poss, good, test, map, indmap, indgood;
    > 
    >      for i in [ 1 .. Length( ComputedPowerMaps( b ) ) ] do
    >        if IsBound( ComputedPowerMaps( b )[i] ) then
    >          PowerMap( s, i );
    >        fi;
    >      od;
    > 
    >      if initmap = fail then
    >        sfusb:= InitFusion( s, b );;
    >      else
    >        sfusb:= initmap;
    >      fi;
    > 
    >      if not TestConsistencyMaps( ComputedPowerMaps( s ), sfusb,
    >                 ComputedPowerMaps( b ) ) then
    >        Error( "inconsistency in power maps!" );
    >      fi;
    > 
    >      poss:= FusionsAllowedByRestrictions( s, b, Irr( s ), ind, sfusb,
    >                 rec( decompose:= true, minamb:= 2, maxamb:= 10^4,
    >                      quick:= false, maxlen:= 10,
    >                      contained:= ContainedPossibleCharacters ) );
    >      indgood:= [];
    >      if ForAll( poss, x -> ForAll( x, IsInt ) ) then
    >        # All candidates in 'poss' are unique.
    >        # Consider only representatives under the symmetry group of 's'.
    >        poss:= RepresentativesFusions( s, poss, Group( () ) );
    > 
    >        # Discard candidates for which the scalar products
    >        # of induced characters are not integral.
    >        good:= [];
    >        test:= ( n  -> IsInt( n ) and 0 <= n );
    >        for map in poss do
    >          indmap:= InducedClassFunctionsByFusionMap( s, b, Irr( s ), map );
    >          if ForAll( indmap,
    >                 x -> ForAll( indmap,
    >                          y -> test( ScalarProduct( b, x, y ) ) ) ) then
    >            Add( good, map );
    >          fi;
    >        od;
    >        poss:= good;
    > 
    >        # Compute those induced characters that arise independent of
    >        # the fusion map.
    >        indgood:= Intersection( List( good,
    >            map -> InducedClassFunctionsByFusionMap( s, b, Irr( s ),
    >                       map ) ) );
    >      fi;
    > 
    >      return rec( maps:= poss, induced:= indgood );
    >    end;;
\end{verbatim}

Our initial characters of $\B$ are the trivial character and the characters
that arise from inducing irreducible characters of cyclic subgroups.

\begin{verbatim}
    gap> knownirr:= [ TrivialCharacter( b ) ];;
    gap> indcyc:= InducedCyclic( b, [ 2 .. NrConjugacyClasses( b ) ], "all" );;
\end{verbatim}

The class fusion from $Fi_{23}$ to $\B$ is determined uniquely
by the available data, and this takes only a few seconds.


\begin{verbatim}
    gap> fi23:= CharacterTable( "Fi23" );;
    gap> fi23fusb:= tryFusion( fi23, b, indcyc, fail );;
    gap> Length( fi23fusb.maps );
    1
    gap> indfi23:= fi23fusb.induced;;
\end{verbatim}

The class fusion from $\Cent_B(2B)$ to $\B$ may be assumed,
see Section~\ref{table_c2b}.

\begin{verbatim}
    gap> b2b:= CharacterTable( "BM2" );;
    gap> b2bfusb:= GetFusionMap( b2b, b );;
    gap> indb2b:= Set( InducedClassFunctionsByFusionMap( b2b, b,
    >                      Irr( b2b ), b2bfusb ) );;
\end{verbatim}

The subgroups $Th$ and $HN.2$ are treated in the same way as $Fi_{23}$.

\begin{verbatim}
    gap> hn2:= CharacterTable( "HN.2" );;
    gap> ind:= Concatenation( indfi23, indb2b, indcyc );;
    gap> hn2fusb:= tryFusion( hn2, b, ind, fail );;
    gap> Length( hn2fusb.maps );
    1
    gap> indhn2:= hn2fusb.induced;;
    gap> th:= CharacterTable( "Th" );;
    gap> ind:= Concatenation( indfi23, indb2b, indcyc );;
    gap> thfusb:= tryFusion( th, b, ind, fail );;
    gap> Length( thfusb.maps );
    1
    gap> indth:= thfusb.induced;;
\end{verbatim}


Now we want to determine the class fusion from $H = 2.{}^2E_6(2).2$ to $\B$.
The approach used above is not feasible in this case.
In order to refine the initial approximation of the class fusion,
we use that $H'$ contains a subgroup of the type $2.Fi_{22}$
that is contained also in a $Fi_{23}$ type subgroup of $\B$.
Note that $H$ is the centralizer of an involution $z$ in $\B$ from the class
\texttt{2A}, and the class \texttt{2A} of $Fi_{23}$ lies in this class.
We may choose our $Fi_{23}$ subgroup such that it contains $z$.
The centralizer of $z$ in $Fi_{23}$ has then the type $2.Fi_{22}$.


Thus we compute the possible class fusions from $2.Fi_{22}$ to $Fi_{23}$
and to $H'$.
The compositions of the former maps with the known fusion from $Fi_{23}$
to $\B$ yields the possible class fusions from $2.Fi_{22}$ to $\B$,
and the compositions of these fusions with the inverses of the latter maps
yield the desired approximations for the fusion from $Fi_{23}$ to $\B$.

\begin{verbatim}
    gap> 2fi22:= CharacterTable( "2.Fi22" );;
    gap> 2fi22fusfi23:= PossibleClassFusions( 2fi22, fi23 );;
    gap> 2fi22fusb:= Set( List( 2fi22fusfi23, map -> CompositionMaps(
    >        fi23fusb.maps[1], map ) ) );;
    gap> Length( 2fi22fusb );
    2
    gap> hh:= CharacterTable( "2.2E6(2)" );;
    gap> 2fi22fushh:= PossibleClassFusions( 2fi22, hh );;
    gap> approxhhfusb:= [];;
    gap> for map1 in 2fi22fushh do
    >      for map2 in 2fi22fusb do
    >        AddSet( approxhhfusb, CompositionMaps( map2, InverseMap( map1 ) ) );
    >      od;
    >    od;
    gap> Length( approxhhfusb );
    4
    gap> inithhfusb:= InitFusion( hh, b );;
    gap> TestConsistencyMaps( ComputedPowerMaps( hh ), inithhfusb,
    >        ComputedPowerMaps( b ) );
    true
    gap> for i in [ 1 .. Length( approxhhfusb ) ] do
    >      if MeetMaps( approxhhfusb[i], inithhfusb ) <> true then
    >        Unbind( approxhhfusb[i] );
    >      fi;
    >    od;
    gap> approxhhfusb:= Compacted( approxhhfusb );;
    gap> Length( approxhhfusb );
    2
\end{verbatim}

We get two initial approximations, and the computation of the
class fusion from $H'$ to $\B$ is now easy.


\begin{verbatim}
    gap> hhfusb:= List( approxhhfusb, map -> tryFusion( hh, b, ind, map ) );;
    gap> List( hhfusb, r -> Length( r.maps ) );
    [ 1, 1 ]
    gap> hhfusb[1] = hhfusb[2];
    true
    gap> hhfusb:= hhfusb[1];;
\end{verbatim}

Thus we have determined the class fusion from $H'$ to $\B$ uniquely.
The next step is to compute the class fusion for the classes in $H$
that do not lie in $H'$.

\begin{verbatim}
    gap> h:= CharacterTable( "2.2E6(2).2" );;
    gap> hhfush:= PossibleClassFusions( hh, h );;
    gap> Length( hhfush );
    4
    gap> approxhfusb:= Set( List( hhfush,
    >        map -> CompositionMaps( hhfusb.maps[1], InverseMap( map ) ) ) );;
    gap> Length( approxhfusb );
    2
    gap> inithfusb:= InitFusion( h, b );;
    gap> TestConsistencyMaps( ComputedPowerMaps( h ), inithfusb,
    >        ComputedPowerMaps( b ) );
    true
    gap> List( approxhfusb, map -> MeetMaps( map, inithfusb ) );
    [ true, true ]
    gap> ind:= Concatenation( indfi23, indb2b, indhn2, hhfusb.induced,
    >              indcyc );;
    gap> hfusb:= List( approxhfusb, map -> tryFusion( h, b, ind, map ) );;
    gap> List( hfusb, r -> Length( r.maps ) );
    [ 1, 1 ]
    gap> hfusb[1].induced = hfusb[2].induced;
    true
    gap> hfusb:= hfusb[1];;
    gap> indh:= hfusb.induced;;
\end{verbatim}


Now we know many induced characters of $\B$.
The $\Z$-lattice that is spanned by these characters contains several
irreducible characters.
Unfortunately, the LLL program in {\GAP} does not find them immediately.

Therefore, we proceed now in two steps.
First,
we assume that $\B$ has a rational ordinary irreducible character $\chi$,
say, of degree $4371$ whose $3$- and $5$-modular restrictions are
the Brauer characters of the representations
which we have used in Section~\ref{invs_B}.
From $\chi$ together with the known induced characters,
we easily compute a list of vectors of norm $1$
such that the input characters are linear combinations of these vectors,
with nonnegative integer coefficients.
In the second step,
we will then \textbf{not} use $\chi$
but we use the vectors found in the first step
as our candidates for the irreducible characters,
which just have to be verified.

Let us start with the first step, and compute the values of $\chi$.
For each representative of order not divisible by $30$,
we compute the Brauer character value from a representation in
characteristic coprime to the element order;
for the remaining classes, we store 'fail' as a preliminary value.

\begin{verbatim}
    gap> chi:= [];;
    gap> for nam in labels do
    >      slp:= SLPForClassName( nam, cycprg, outputnames );
    >      ord:= Int( Filtered( nam, IsDigitChar ) );
    >      if ord mod 2 <> 0 then
    >        val:= 1 + BrauerCharacterValue(
    >                      ResultOfStraightLineProgram( slp, gens_2 ) );
    >      elif ord mod 3 <> 0 then
    >        val:= BrauerCharacterValue(
    >                  ResultOfStraightLineProgram( slp, gens_3 ) );
    >      elif ord mod 5 <> 0 then
    >        val:= BrauerCharacterValue(
    >                  ResultOfStraightLineProgram( slp, gens_5 ) );
    >      else
    >        val:= fail;
    >      fi;
    >      Add( chi, val );
    >    od;
    gap> chi;
    [ 4371, -493, 275, 19, -53, 78, -3, 51, 3, -1, -5, -77, 11, -13, 19, -21, 35, 
      -4, 21, 20, 4, 1, 14, -7, -8, 13, -2, 1, 5, -34, 10, -5, 3, -1, -9, 11, 7, 
      -5, 3, 7, -1, -1, -1, -1, -21, -3, 6, -3, -1, 5, 0, 7, 4, 4, 2, -3, 0, -3, 
      2, -4, 5, -12, 0, 1, 1, -1, -1, -2, -3, 4, 5, -6, 6, 3, 3, -2, -4, 4, 2, 
      -10, 2, 3, -5, -1, 3, 1, -1, 3, 1, 2, 2, 2, -2, 1, 2, -2, 1, -1, 1, 0, 0, 
      1, -1, 2, -2, 3, -7, 1, 0, 2, 1, 1, 0, -3, 0, 1, -4, 2, -2, -2, 2, -1, 1, 
      0, -1, 1, 1, -1, 0, 2, -2, 0, 0, 0, fail, fail, fail, fail, fail, fail, 
      fail, 0, 1, -1, 1, -2, 0, 0, -1, 0, 0, 1, 0, -1, 0, 1, -1, -3, -1, -1, 1, 
      0, -1, 0, 0, -2, -1, -1, 0, fail, fail, fail, -1, 0 ]
\end{verbatim}


Next we transfer these values to the class positions in the character table.
For that, we compute the mapping from the character table head to the labels.

\begin{verbatim}
    gap> bfuslabels:= [];;
    gap> for i in [ 1 .. Length( libBnames ) ] do
    >      poss:= Filtered( labels, x -> Filtered( libBnames[i], IsDigitChar )
    >             = Filtered( x, IsDigitChar ) and
    >             ForAny( Filtered( x, IsAlphaChar ), y -> y in libBnames[i] ) );
    >      if Length( poss ) <> 1 then
    >        Print( "problem for ", libBnames[i], "\n" );
    >      fi;
    >      bfuslabels[i]:= Position( labels, poss[1] );
    >    od;
\end{verbatim}

The character values are unknown for eight classes of element order $30$
and three classes of element order $60$.

\begin{verbatim}
    gap> chi:= chi{ bfuslabels };
    [ 4371, -493, 275, -53, 19, 78, -3, -77, 51, 19, -21, 35, -13, 11, 3, -1, -5, 
      21, -4, -34, 20, 14, -7, 4, -8, 5, -2, 13, 1, 1, 10, -21, 7, -9, 11, -1, 
      -5, -5, 3, -1, 3, 7, -1, -1, -1, -3, 6, 7, 5, -3, 0, -1, 4, 4, -12, -3, 4, 
      6, -6, 5, -2, 0, -4, 2, -3, 2, 1, -1, 5, 0, -3, 1, 3, -1, 3, -10, 4, -4, 2, 
      -2, 3, 2, 3, -5, -1, 3, -1, 3, 1, 1, 2, 2, 2, 2, -2, 1, -2, 1, -7, -2, 3, 
      1, -1, 1, 0, -1, 2, 0, 1, 2, 0, 1, 1, -2, 0, 2, -4, 0, -3, 2, 1, -2, 0, 1, 
      1, -1, -1, 1, -1, 1, 0, 2, -2, 0, 0, 0, fail, fail, fail, fail, fail, fail, 
      fail, fail, 0, 0, 1, 1, -1, -1, 1, -2, 0, 0, 0, 0, 0, -1, 1, 0, -3, 1, -1, 
      -1, 0, -1, 1, -1, 0, -1, -1, 0, 0, 0, -2, -1, -1, 0, 0, fail, fail, fail, 
      -1, 0 ]
    gap> failpos:= Positions( last, fail );
    [ 137, 138, 139, 140, 141, 142, 143, 144, 180, 181, 182 ]
    gap> OrdersClassRepresentatives( b ){ failpos };
    [ 30, 30, 30, 30, 30, 30, 30, 30, 60, 60, 60 ]
\end{verbatim}

In order to compute the missing values,
we use the same idea as in Section~\ref{table_c2b}.
We know that these values are integers.
For all classes in question except one,
the class lengths are at least $|B| / 360$,
thus the absolute value of the character value at these classes
cannot exceed $5$.

\begin{verbatim}
    gap> SizesCentralizers( b ){ failpos };
    [ 3600, 360, 360, 240, 240, 120, 60, 60, 120, 120, 60 ]
    gap> bclasses:= SizesConjugacyClasses( b );;
    gap> normpart:= Sum(
    >        List( Difference( [ 1 .. Length( bclasses ) ], failpos ),
    >              i -> bclasses[i] * chi[i]^2 ) );;
    gap> normpart + 6^2 * Size( b ) / 360 > Size( b );
    true
\end{verbatim}

Since $\chi( g^p ) \equiv \chi(g) \pmod{p}$ for $p \in \{ 3, 5 \}$
and since the values $\chi( g^p )$ are known,
we know the congruence class of $\chi( g )$ modulo $15$ for all missing
classes except one.
This determines $\chi$.

\begin{verbatim}
    gap> for i in failpos do
    >      v:= ChineseRem( [ 3, 5 ], [ chi[ PowerMap( b, 3, i ) ],
    >                                  chi[ PowerMap( b, 5, i ) ] ] );
    >      if v > 5 then
    >        v:= v - 15;
    >      fi;
    >      chi[i]:= v;
    >    od;
    gap> chi{ failpos };
    [ -5, 1, -3, -1, -1, -2, 0, 0, -1, -1, 0 ]
    gap> Sum( List( [ 1 .. Length( bclasses ) ],
    >               i -> bclasses[i] * chi[i]^2 ) ) = Size( b );
    true
\end{verbatim}

Now we compute candidates for the irreducibles,
under the assumption that $\chi$ is a character.
The following simple function is later used in a loop.
It takes a character table \texttt{tbl} and three lists of (virtual)
characters of this table:
\texttt{knownirr} contains known irreducibles,
\texttt{newirr} contains newly found irreducibles,
and \texttt{knownvirt} contains virtual characters.
The idea is as follows.
The list \texttt{knownirr} gets extended (in place) by \texttt{newirr},
symmetrizations of the characters in \texttt{newirr} and tensor products
of the characters in \texttt{knownirr} with those in \texttt{newirr} are
created, and the concatenation of \texttt{knownvirt} and these characters
gets reduced with \texttt{knownirr} and \texttt{newirr};
this process is iterated as long as new irreducible characters are found
in the reduction step, with \texttt{newirr} replaced by these characters;
the function returns the list of non-irreducible remainders of
\texttt{knownvirt}.

\begin{verbatim}
    gap> ExtendCandidates:= function( tbl, knownirr, newirr, knownvirt )
    >      while 0 < Length( newirr ) do
    >        Append( knownirr, newirr );
    >        knownvirt:= Concatenation( knownvirt,
    >                        Concatenation( List( [ 2, 3, 5 ],
    >                            p -> Symmetrizations( b, newirr, p ) ) ),
    >                        Set( Tensored( knownirr, newirr ) ) );
    >        knownvirt:= Reduced( tbl, knownirr, knownvirt );
    >        newirr:= knownvirt.irreducibles;
    >        knownvirt:= knownvirt.remainders;
    >      od;
    > 
    >      Print( "#I  ", Length( knownirr ), " irreducibles known\n" );
    >      return knownvirt;
    >    end;;
\end{verbatim}

We start with the trivial character of $\B$ and $\chi$,
and note that the antisymmetric square of $\chi$ is irreducible.
We initialize \texttt{knownvirt} with the union of the induced characters
which were generated above.

\begin{verbatim}
    gap> psi:= AntiSymmetricParts( b, [ chi ], 2 )[1];;
    gap> ScalarProduct( b, psi, psi );
    1
    gap> testind:= Concatenation( indfi23, indb2b, indhn2, indh, indth, indcyc );;
\end{verbatim}

Reducing the characters with our three irreducible characters
and applying the LLL algorithm to the remainders yields $7$
new irreducibles.

\begin{verbatim}
    gap> knownirr:= [];;
    gap> l:= ExtendCandidates( b, knownirr,
    >         [ TrivialCharacter( b ), chi, psi ], testind );;
    #I  3 irreducibles known
    gap> lll:= LLL( b, l, 99/100 );;  Length( lll.irreducibles );
    7
\end{verbatim}

We are lucky,
a short loop of reductions with the new irreducibles and LLL reduction
yields the complete list of irreducible characters.
The ch

\begin{verbatim}
    gap> n:= NrConjugacyClasses( b );;
    gap> while 0 < Length( lll.irreducibles ) and
    >       Length( lll.irreducibles ) + Length( knownirr ) < n do
    >      l:= ExtendCandidates( b, knownirr,
    >                            lll.irreducibles, lll.remainders );;
    >      lll:= LLL( b, l, 99/100 );;
    >    od;
    #I  11 irreducibles known
    #I  17 irreducibles known
    gap> Append( knownirr, lll.irreducibles );
\end{verbatim}

The irreducible characters found this way coincide with the ones from the
{\ATLAS} table of $\B$.

\begin{verbatim}
    gap> Set( irr_atlas ) = Set( knownirr );
    true
\end{verbatim}

Now comes step two.
As stated above, all what remains is to verify the candidate vectors,
where we are allowed to use the induced characters but not $\chi$.

For each candidate vector, we compute whether it occurs
in the $\Z$-span of the induced characters,
and if yes, we compute the coefficients of a linear combination in terms
of them.
This way, we find/verify the first $30$ irreducible characters of $\B$.




\begin{verbatim}
    gap> ind:= Concatenation( indfi23, indb2b, indhn2, indh, indcyc );;
    gap> mat:= MatScalarProducts( b, knownirr, ind );;
    gap> irr:= [ TrivialCharacter( b ) ];;
    gap> one:= IdentityMat( NrConjugacyClasses( b ) );;
    gap> for i in [ 2 .. Length( one ) ] do
    >      coeffs:= SolutionIntMat( mat, one[i] );
    >      if coeffs <> fail and ForAll( coeffs, IsInt ) then
    >        Add( irr, coeffs * ind );
    >      fi;
    >    od;
    gap> Length( irr );
    30
    gap> Set( List( irr, chi -> ScalarProduct( b, chi, chi ) ) );
    [ 1 ]
    gap> ForAll( irr, chi -> chi[1] > 0 );
    true
\end{verbatim}

In order to get the missing irreducible characters of $\B$,
we add symmetrizations and tensor products of the known irreducibles
to the list of induced characters, and project onto the orthogonal space
of the space that is spanned by the known irreducibles.
No new irreducible characters are found directly this way,
but the \verb|`oracle|' that was used above tells us how to express the missing
irreducibles as integral linear combinations of the known characters.

\begin{verbatim}
    gap> sym:= Symmetrizations( b, irr, 2 );;
    gap> Append( sym, Symmetrizations( b, irr, 3 ) );
    gap> Append( sym, Symmetrizations( b, irr, 5 ) );
    gap> ten:= Set( Tensored( irr, irr ) );;
    gap> cand:= Reduced( b, irr, Concatenation( sym, ten, ind ) );;
    gap> cand.irreducibles;
    [  ]
    gap> cand:= cand.remainders;;
    gap> newirr:= [];;
    gap> mat:= MatScalarProducts( b, knownirr, cand );;
    gap> for i in [ 2 .. Length( one ) ] do
    >      coeffs:= SolutionIntMat( mat, one[i] );
    >      if coeffs <> fail and ForAll( coeffs, IsInt ) then
    >        Add( newirr, coeffs * cand );
    >      fi;
    >    od;
    gap> Length( newirr );
    154
    gap> Set( List( newirr, chi -> ScalarProduct( b, chi, chi ) ) );
    [ 1 ]
    gap> ForAll( newirr, chi -> chi[1] > 0 );
    true
\end{verbatim}

It is not surprising that
the irreducible characters found this way coincide with the ones from the
{\ATLAS} table of $\B$.

\begin{verbatim}
    gap> Append( irr, newirr );
    gap> Set( irr_atlas ) = Set( irr );
    true
\end{verbatim}

\section{Appendix: Standardizing the generators of $Co_2$}\label{slp_co2}

In Section~\ref{table_c2b}, we have restricted the $3$-modular
$4371$-dimensional representation of $\B$ to a subgroup $C$ of the structure
$2^{1+22}.Co_2$, and obtained a $23$-dimensional composition factor
with a faithful action of the factor group $Co_2$,
with generators $x$ and $y$, say.
Our aim is to find short words in terms of $x$ and $y$
that yield \emph{standard} generators for $Co_2$,
that is, elements $a$, $b$ from the classes \texttt{2A} and \texttt{5A}
of $Co_2$, with the properties that the product $a b$ has order $28$
and that $a$ and $b$ generate the full group $Co_2$.

The classes \texttt{2A} and \texttt{5A} of $Co_2$ are determined by the
Brauer character values $-9$ and $-2$, respectively,
in the unique irreducible $23$-dimensional $3$-modular representation
of $Co_2$.

\begin{verbatim}
    gap> t:= CharacterTable( "Co2" ) mod 3;;
    gap> dim23:= Filtered( Irr( t ), chi -> chi[1] = 23 );;
    gap> Length( dim23 );
    1
    gap> pos:= PositionsProperty( OrdersClassRepresentatives( t ),
    >                             x -> x in [ 2, 5 ] );
    [ 2, 3, 4, 12, 13 ]
    gap> dim23[1]{ pos };
    [ -9, 7, -1, -2, 3 ]
\end{verbatim}

The $12$-th power of the second generator yields a \texttt{2A} element.

\begin{verbatim}
    gap> List( co2gens, Order );
    [ 2, 24 ]
    gap> BrauerCharacterValue( co2gens[1] );
    7
    gap> a:= co2gens[2]^12;;
    gap> BrauerCharacterValue( a );
    -9
\end{verbatim}

A short word that defines a \texttt{5A} element is $(y^4 x)^4$.
It can be found as follows.

\begin{verbatim}
    gap> f:= FreeMonoid( 2 );;
    gap> fgens:= GeneratorsOfMonoid( f );
    [ m1, m2 ]
    gap> for w in Iterator( f ) do
    >      m:= MappedWord( w, fgens, co2gens );
    >      ord:= Order( m );
    >      if ord mod 5 = 0 then
    >        cand:= m^( ord / 5 );
    >        if BrauerCharacterValue( cand ) = -2 then
    >          break;
    >        fi;
    >      fi;
    >    od;
    gap> w;
    m2^4*m1
    gap> ord;
    20
\end{verbatim}

Similarly, we find a conjugate of the \texttt{5A} element such that
the product has order $28$.

\begin{verbatim}
    gap> for w in Iterator( f ) do
    >      m:= MappedWord( w, fgens, co2gens );
    >      b:= cand ^ m;
    >      if Order( a * b ) = 28 then
    >        break;
    >      fi;
    >    od;
    gap> w;
    m2*(m2*m1)^3
\end{verbatim}

In order to show that the two elements \texttt{a} and \texttt{b}
really generate $Co_2$,
we use the fact that no proper subgroup of $Co_2$ contains elements
of the orders $23$ and $30$.

\begin{verbatim}
    gap> t:= CharacterTable( "Co2" );;
    gap> mx:= List( Maxes( t ), CharacterTable );;
    gap> ForAny( List( Maxes( t ), CharacterTable ),
    >            x -> IsSubset( OrdersClassRepresentatives( x ),
    >                           [ 23, 30 ] ) );
    false
\end{verbatim}

It suffices to find products of the generators that have these orders.
We compute random elements until we are successful.

\begin{verbatim}
    gap> u:= Group( a, b );;
    gap> found:= [];;
    gap> repeat
    >      x:= Order( PseudoRandom( u ) );
    >      if x in [ 23, 30 ] then
    >        AddSet( found, x );
    >      fi;
    >    until Length( found ) = 2;
    gap> found;
    [ 23, 30 ]
\end{verbatim}


\section{Appendix: Words for generators of the kernel $2^{22}$}%
\label{slp_kernel}

We assume (see Section~\ref{table_c2b})
that we have permutation generators \texttt{stdperms} on $4060$ points
for the group $P \cong 2^{22}.Co_2$ such that mapping them to
standard generators of $Co_2$ defines an epimorphism.
Our aim is to find (short) words in terms of \texttt{stdperms} for
generators of the elementary abelian normal subgroup of order $2^{22}$.

For that, we work in parallel with standard generators \texttt{co2permgens}
of $Co_2$.

\begin{verbatim}
    gap> co2permgens:= GeneratorsOfGroup( AtlasGroup( "Co2" ) );;
    gap> f:= FreeMonoid( 2 );;
    gap> fgens:= GeneratorsOfMonoid( f );;
    gap> kernelgens:= [];;
    gap> kernelwords:= [];;
    gap> kernel:= Group( () );;
    gap> for w in Iterator( f ) do
    >      m1:= MappedWord( w, fgens, stdperms );
    >      m2:= MappedWord( w, fgens, co2permgens );
    >      ord1:= Order( m1 );
    >      ord2:= Order( m2 );
    >      if ord1 <> ord2 then
    >        cand:= m1 ^ ord2;
    >        if not cand in kernel then
    >          Add( kernelgens, cand );
    >          Add( kernelwords, [ w, ord2 ] );
    >          kernel:= ClosureGroup( kernel, cand );
    >        fi;
    >      fi;
    >      if Length( kernelgens ) = 22 then
    >        break;
    >      fi;
    >    od;
\end{verbatim}

The kernel generators obtained with the above procedure are equal to
the results of the straight line program \texttt{slp\_ker}
from Section~\ref{table_c2b}.

\begin{verbatim}
    gap> ResultOfStraightLineProgram( slp_ker, fgens )
    >    = List( kernelwords, pair -> pair[1]^pair[2] );
    true
\end{verbatim}

\section{Appendix: Words for class representatives of $2^{22}.Co_2$}%
\label{slp_classreps}

The group $P \cong 2^{22}.Co_2$ from Section~\ref{table_c2b}
has $388$ conjugacy classes, as we can compute from the permutation
representation on $4060$ points.
Our aim is to find (short) words in terms of the generators \texttt{stdperms}
for conjugacy class representatives of $P$.

For that, we first compute preimages in $P$ of class representatives
of its factor group $Co_2$,
using the straight line program \texttt{slp\_co2classreps} from~\cite{AGRv3},
and use them to initialize lists of known class representatives and of words
for each class of $Co_2$.
In the lists of words, we record the indices of kernel generators with
which we have to multiply the preimage; thus the preimages themselves
are denoted by empty lists.

\begin{verbatim}
    gap> factccls:= ResultOfStraightLineProgram( slp_co2classreps.program,
    >                   stdperms );;
    gap> classreps:= List( factccls, x -> [ x ] );;
    gap> classwords:= List( factccls, x -> [ [] ] );;
\end{verbatim}

The preimage of the identity element in $Co_2$ is not the identity in $P$.
Since it would be hard to get the identity as a product of this preimage
with a product of kernel generators,
we add the identity element by hand, and denote it by the word \texttt{[ 0 ]}.

\begin{verbatim}
    gap> classreps[1]:= Concatenation( [ () ], classreps[1] );;
    gap> classwords[1]:= Concatenation( [ [ 0 ] ], classwords[1] );;
\end{verbatim}

Now we multiply the class representatives by words (of increasing length)
in terms of the kernel generators
and add the elements if they yields new classes,
until we have found enough class representatives.

The conjugacy checks in the following piece of {\GAP} code
cannot be executed with {\GAP}~4.9.3  or earlier versions in reasonable time.
We have done these computations via delegations to {\MAGMA}~\cite{Magma},
using the auxiliary function \texttt{IsConjugateViaMagma}.
(We hope that eventually the necessary functionality will become available
also in {\GAP}.)

\begin{verbatim}
    gap> if CTblLib.IsMagmaAvailable() then
    >      IsConjugateViaMagma:= function( permgroup, pi, known )
    >        local path, inputs, str, out, result, pos;
    > 
    >        path:= UserPreference( "CTblLib", "MagmaPath" );
    >        inputs:= [ CTblLib.MagmaStringOfPermGroup( permgroup, "G" ),
    >                   Concatenation( "p:= G!", String( pi ), ";" ),
    >                   "l:= [" ];
    >        if Length( known ) > 0 then
    >          Append( inputs,
    >                List( known, p -> Concatenation( "G!", String( p ), "," ) ) );
    >          Remove( inputs[ Length( inputs ) ] );
    >        fi;
    >        Append( inputs,
    >                [ "];",
    >                  "conj:= false;",
    >                  "for q in l do",
    >                  "  conj:= IsConjugate( G, p, q );",
    >                  "  if conj then break; end if;",
    >                  "end for;",
    >                  "if conj then",
    >                  "  print true;",
    >                  "else",       # Do not call '# Class( G, p );'.
    >                  "  print \"#\", # Centralizer( G, p );",
    >                  "end if;" ] );
    >        str:= "";
    >        out:= OutputTextString( str, true );
    >        result:= Process( DirectoryCurrent(), path,
    >            InputTextString( JoinStringsWithSeparator( inputs, "\n" ) ),
    >            out, [] );
    >        CloseStream( out );
    >        if result <> 0 then
    >          Error( "Process returned ", result );
    >        fi;
    >        pos:= PositionSublist( str, "\n# " );
    >        if pos <> fail then
    >          return Int( str{ [ pos + 3 .. Position( str, '\n', pos+1 )-1 ] } );
    >        elif PositionSublist( str, "true" ) <> fail then
    >          # 'pi' is conjugate to a perm. in 'known'
    >          return true;
    >        else
    >          Error( "Magma failed" );
    >        fi;
    >      end;
    >      co2tbl:= CharacterTable( "Co2" );;
    >      g:= Group( stdperms );;
    >      for i in [ 1 .. Length( factccls ) ] do
    >        iclasses:= List( classreps[i], x -> ConjugacyClass( g, x ) );
    >        sum:= 2^22 * SizesConjugacyClasses( co2tbl )[i]
    >              - Sum( List( iclasses, Size ) );
    >        len:= 1;
    >        while sum > 0 do
    >          for tup in Combinations( [ 1 .. 22 ], len ) do
    >            cand:= factccls[i] * Product( kernelgens{ tup } );
    >            # We call Magma anyhow in order to compute the class length.
    >            conj:= IsConjugateViaMagma( g, cand, classreps[i] );
    >            if conj <> true then
    >              Add( classreps[i], cand );
    >              Add( classwords[i], tup );
    >              sum:= sum - Size( g ) / conj;
    >              if sum <= 0 then
    >                break;
    >              fi;
    >            fi;
    >          od;
    >          len:= len + 1;;
    >        od;
    >      od;
    >    fi;
\end{verbatim}

Let us check whether the class representatives fit to the ones
used in Section~\ref{table_c2b} (up to ordering).

\begin{verbatim}
    gap> if CTblLib.IsMagmaAvailable() then
    >      Print( ForAll( [ 1 .. Length( factccls ) ],
    >        i -> Set( classrepsinfo[i][2] ) = Set( classwords[i] ) ), "\n" );
    true
    >    else
    >      Print( "Magma not available, no check of class representatives\n" );
    >    fi;
\end{verbatim}


\section{Appendix: About the character table of $2^{9+16}.S_8(2)$}%
\label{2Dnormalizer}

As has been stated in Section~\ref{sect:classes},
we use the character table of a maximal subgroup of $\B$
that normalizes an elementary abelian group of order $2^8$
instead of the table of the \texttt{2D} normalizer in $\B$.
The character table of this maximal subgroup,
of the structure $2^{9+16}.S_8(2)$,
had been computed from a matrix representation of the group.
We verify that this matrix group can indeed be obtained from the
restriction of one of our certified matrix representations of $\B$.

First we restrict our representation of $\B$ over the field with
two elements to the fourth maximal subgroup.

\begin{verbatim}
    gap> prg:= AtlasProgram( "B", "maxes", 4 );;
    gap> gens:= ResultOfStraightLineProgram( prg.program, gens_2 );;
\end{verbatim}

Next we compute a $180$ dimensional representation of this group.
(The fact that this representation is faithful follows from the
computations with it; this is beyond the scope of this note.)

\begin{verbatim}
    gap> m:= GModuleByMats( gens, GF(2) );;
    gap> a:= gens[1];;  b:= gens[2];;
    gap> mat:= a^2 + a*b^2 + b^3 * a * b^2 + b*a + b^3 * a;;
    gap> nsp:= NullspaceMat( mat );;  Length( nsp );
    5
    gap> s:= MTX.SubGModule( m, nsp[4] );;
    gap> ind:= MTX.InducedActionSubmodule( m, s );;
    gap> MTX.Dimension( ind );
    206
    gap> css:= MTX.BasesCompositionSeries( ind );;
    gap> List( css, Length );
    [ 0, 26, 27, 35, 163, 171, 172, 198, 206 ]
    gap> ind2:= MTX.InducedActionFactorModule( ind, css[2] );;
    gap> MTX.Dimension( ind2 );
    180
\end{verbatim}

The word used above has a $3$ dimensional nullspace
on the $180$ dimensional module.
We try each nonzero vector in that space as a seed vector,
and consider the standard bases of the submodules generated by them;
for that, we apply the function \texttt{StdBasis} from Section~\ref{pres_B}.

\begin{verbatim}
    gap> gensnew:= ind2.generators;;
    gap> a:= gensnew[1];;  b:= gensnew[2];;
    gap> mat:= a^2 + a*b^2 + b^3 * a * b^2 + b*a + b^3 * a;;
    gap> nsp:= NullspaceMat( mat );;  Length( nsp );
    3
    gap> stdbasnew:= List( NormedRowVectors( VectorSpace( GF(2), nsp ) ),
    >                      seed -> StdBasis( GF(2), gensnew, seed ) );;
    gap> List( stdbasnew, Length );
    [ 180, 180, 9, 180, 145, 145, 180 ]
    gap> stdbasnew:= Filtered( stdbasnew, x -> Length( x ) = 180 );;
    gap> stdgensnew:= List( stdbasnew, x -> List( gensnew, m -> x*m*x^-1 ) );;
\end{verbatim}

We repeat the same process with the matrix generators from which
the character table of the maximal subgroup $2^{9+16}.S_8(2)$ of $\B$
has been computed.


\begin{verbatim}
    gap> gensold:= OneAtlasGeneratingSet( "2^(9+16).S8(2)", Dimension, 180 );;
    gap> gensold.repname;
    "Bmax4G0-f2r180B0"
    gap> gensold:= gensold.generators;;
    gap> a:= gensold[1];;  b:= gensold[2];;
    gap> mat:= a^2 + a*b^2 + b^3 * a * b^2 + b*a + b^3 * a;;
    gap> nsp:= NullspaceMat( mat );;  Length( nsp );
    3
    gap> stdbasold:= List( NormedRowVectors( VectorSpace( GF(2), nsp ) ),
    >                      seed -> StdBasis( GF(2), gensold, seed ) );;
    gap> List( stdbasold, Length) ;
    [ 180, 180, 9, 180, 145, 145, 180 ]
    gap> stdbasold:= Filtered( stdbasold, x -> Length( x ) = 180 );;
    gap> stdgensold:= List( stdbasold, x -> List( gensold, m -> x*m*x^-1 ) );;
\end{verbatim}

It turns out that the normal forms obtained this way coincide.

\begin{verbatim}
    gap> List( stdgensnew, x -> Position( stdgensold, x ) );
    [ 1, 2, 4, 3 ]
\end{verbatim}


\newcommand{\etalchar}[1]{$^{#1}$}
\providecommand{\bysame}{\leavevmode\hbox to3em{\hrulefill}\thinspace}
\providecommand{\MR}{\relax\ifhmode\unskip\space\fi MR }
\providecommand{\MRhref}[2]{%
  \href{http://www.ams.org/mathscinet-getitem?mr=#1}{#2}
}
\providecommand{\href}[2]{#2}


\end{document}